\documentclass[preprint,12pt]{elsarticle}

\usepackage{amssymb}
\usepackage{amsmath}
\usepackage{hyperref}
\usepackage{lineno}

\usepackage{bm} 
\usepackage{mathtools}
\usepackage{graphicx}
\usepackage{subcaption}
\graphicspath{{figs/}} 
\hypersetup{pdfauthor=author}
\usepackage{mismath}

\usepackage{tikz}
\journal{IMA Journal of Numerical Analysis}
\newcommand{\dx}{\,\mathrm{d}}
\newcommand{\Real}{\mathbb{R}}
\newcommand{\Complex}{\mathbb{C}}

\DeclarePairedDelimiter{\RoundBrackets}{(}{)}
\DeclarePairedDelimiter{\CurlyBrackets}{\{}{\}}

\usepackage{stmaryrd}

\DeclarePairedDelimiter{\Norm}{\lVert}{\rVert}

\usepackage{amsthm}
\usepackage{cleveref}
\newtheorem{theorem}{Theorem}[section]

\newtheorem{lemma}[theorem]{Lemma}
\newtheorem*{remark}{Remark}

\usepackage{array}
\newcolumntype{L}[1]{>{\raggedright\let\newline\\\arraybackslash\hspace{0pt}}m{#1}}
\newcolumntype{C}[1]{>{\centering\let\newline\\\arraybackslash\hspace{0pt}}m{#1}}
\newcolumntype{R}[1]{>{\raggedleft\let\newline\\\arraybackslash\hspace{0pt}}m{#1}}

\usepackage{algorithm}
\usepackage{algcompatible}
\begin{document}

\begin{frontmatter}


\author[a]{Xingguang Jin\corref{cor1}}
\ead{xgjin@math.cuhk.edu.hk}
\cortext[cor1]{Corresponding author.}
\affiliation[a]{organization={Department of Mathematics, The Chinese University of Hong Kong}, addressline={Shatin}, country={Hong Kong SAR}.}
\author[a]{Changqing Ye}\author[a]{Eric T. Chung}
\title{Robust Multiscale Methods for Helmholtz equations in high contrast heterogeneous media} 

\begin{abstract}
In this paper, we provide the constraint energy minimization generalized multiscale finite element method (CEM-GMsFEM) to solve Helmholtz equations in heterogeneous medium. This novel multiscale method is specifically designed to overcome problems related to pollution effect, high-contrast coefficients, and the loss of hermiticity of operators. We establish the inf-sup stability and give an a priori error estimate for this method under a number of established assumptions and resolution conditions. The theoretical results are validated by a set of numerical tests, which further show that the multiscale technique can effectively capture pertinent physical phenomena.
\end{abstract}


\end{frontmatter}




\section{Introduction}
\linenumbers
The Helmholtz equation models wave propagation and scattering phenomena in the frequency domain arising in a variety of science and engineering applications, including seismic imaging that helps to model and analyze how seismic waves propagate through the Earth's subsurface, medical ultrasound technologies to simulate and understand the behavior of ultrasound waves in human tissues for diagnostic imaging purposes, and underwater acoustics to study the propagation of sound waves in the ocean and predict their behavior for sonar applications, underwater communication, and marine research \cite{Romanow2008, Lahaye2017}. 
Designing robust and accurate numerical methods for solving the Helmholtz equation can still be challenging, particularly when dealing with high-wavenumber problems or heterogeneous coefficients \cite{Babu1997}. It is well-known that Galerkin Finite element Method (FEM) leads to quasi-optimal error estimates with respect to the degrees of freedoms which is affected by the wavenumber when applied to the Helmholtz equations. 
As a result, the mesh size should be chosen small enough and also the arising system of linear equations is highly indefinite such that the solution process becomes too expensive for the large wavenumber \cite{Ihlenburg1997}. 
The properties of material are usually assumed to have constant density and same speed of sound when study acoustic wave propagation, and for the Helmholtz equations with constant coefficients or low contrast coefficients, pollution effect can be resolved by  many well developed and designed approaches such as the partition of unity finite element methods \cite{Luis2020}, the least squares finite element methods \cite{CHANG19901}, the generalized finite element methods \cite{John2006}, the hybridized discontinuous Galerkin methods \cite{Chen2013}, the interior penalty discontinuous Galerkin methods \cite{Feng2009} and $hp$-FEM \cite{MELENK2010, Melenk2011}. When dealing with complex materials, such as composites, where multiple scales are present within the domain, the challenges become even more significant. In such cases, multiscale methods can be employed to resolve the micro scales and capture the effective behavior of the material. 
These methods aim to bridge the gap between the fine-scale details and the macroscopic behavior by incorporating appropriate scale decomposition techniques and coupling strategies. 

Multiscale methods have long been developed to solve the difficulties associated with the rough coefficients in elliptic equations \cite{Eric2023,FU2019,maalqvist2014}. A commonly shared idea of these methods is to encode fine-scale information into the basis functions of finite element methods. Once the fine-scale information is encoded into the basis functions, the original problem can be solved using multiscale finite element spaces. These spaces have reduced dimensions compared to the default FEM spaces. The reduction in dimensions is achieved by effectively representing the solution behavior at multiple scales using a smaller set of basis functions. These methods have significantly advanced the field of multiscale computational methods and provided efficient and accurate tools for solving problems with multiscale characteristics. These methods have shown promise in capturing the behavior of solutions at different scales and providing accurate approximations in problems with high contrast and heterogeneity. These methods have also been applied to solve Helmholtz equations in heterogeneous domains efficiently, such as Localized Orthogonal Decomposition (LOD) \cite{Peterseim2019, Peterseim2014} method, Super-LOD method \cite{Philip2021},  Multiscale Petrov--Galerkin Method \cite{brown2017}, Generalized Multiscale Finite Element methods (GMsFEM) \cite{Gavrilieva2019}, Edge Multiscale Interior Penalty Discontinuous Galerkin method \cite{Fu2021}, Multiscale spectral generalized finite element method \cite{Ma2023}, Exponentially convergent multiscale methods \cite{chen2023} etc. Recently, a novel multiscale method named Constraint Energy Minimization Generalized Multiscale Finite Element Method (CEM-GMsFEM) is initially developed by \cite{Chung2018} which is aimed for the high-contrast problems  and it has been successfully applied to various partial differential equations arising from practical applications, see, e.g., \cite{Zhao2023,WANG2024,WANG2021}. This approach and the LOD method have certain similarities. For instance, they both rely on the exponential decay features of basis functions and require mesh sizes that are dependent on the oversampling regions in order to achieve the appropriate convergence rate \cite{Henning2021}. Instead of using quasi-interpolation operators to split the solution into macroscopic and microscopic components in LOD, we employ element-wise projections in the implementation of CEM-GMsFEM method, which is the core idea from GMsFEM \cite{Yalchin2013}. To the best of our knowledge, CEM-GMsFEM has not been widely explored to solve Heterogeneous Helmholtz equations. In this paper, we focus on the specific setting of the diffusion coefficients which take either $\varepsilon^2$ on one part of the domain or a small value on the complement of the domain due to the fact that the effects observed in time-harmonic wave propagation can differ depending on whether the diffusion coefficient is small or not \cite{John2008, Peterseim2019}. In the original CEM-GMsFEM, the newly obtained trial and test spaces of multiscale basis functions are the same; however, because the Helmholtz equations are inherently indefinite, we take a cue from the Multiscale Petrov-Galerkin method \cite{Maier2020,brown2017} and form distinct trial and test multiscale spaces by leveraging standard procedures of CEM-GMsFEM in order to  develop a tailored approach for solving the Helmholtz equation in such scenarios.

The aim of this paper is to explore an efficient method to solve Heterogeneous Helmholtz problems by using the novel multiscale model reduction skills coming from CEM-GMsFEM, which beyond the need of periodic coefficients  \cite{brown2017} or other requirements of the coefficients structures. In the analysis section, we establish a resolution condition and build the inf-sup condition of both global problem and multiscale problem to secure their well-posedness. Subsequently, the exponential decay properties of basis functions are demonstrated, and ultimately, we obtain the error estimate of our multiscale method with the desired the convergence rate. For the first time, we present the evidence supporting of the convergence of CEM-GMsFEM  for the Helmholtz equations in heterogeneous media. The numerical simulation section displays the three experiments correspond to three kinds of media, which supports the effectiveness of CEM-GMsFEM and the pollution effect is resolved by using the coarser mesh size to achieve the quasi-optimal convergence. We evaluate the relative error of the CEM-GMsFEM method with respect to different coarse mesh sizes and different oversampling layers in \cref{tab:3,tab:4,tab:5}. The oversampling layers refer to additional layers of elements or degrees of freedom surrounding the coarse mesh, which capture the fine-scale details and improve the accuracy of the approximation. The results of the experiments indicate that the relative errors are influenced by the choice of oversampling layers, which distinguishes the CEM-GMsFEM method from traditional FEM. This suggests that the oversampling layers play a significant role in capturing the fine-scale information and reducing the approximation error. We also compare the relative errors obtained with different oversampling layer configurations to demonstrate the impact of these layers on the accuracy of the method.

By highlighting the influence of oversampling layers on the relative errors, we emphasize the advantage of the CEM-GMsFEM method over traditional FEM in capturing fine-scale information and improving the accuracy of the solution. This finding further supports the effectiveness and efficiency of the proposed method in solving Helmholtz equations in heterogeneous domains.

Our paper is organized as follows. In \cref{sec:pre}, we present the weak formulations of  Helmholtz equations with homogeneous Robin boundary conditions as well as introducing the notations of grids. Then in \cref{sec:cemgms} we construct  multiscale basis functions to construct  multiscale trial space and multiscale test space of our CEM-GMsFEM respectively. In \cref{sec:anal}, we analyze the convergence of CEM-GMsFEM and derive the error estimates. In \cref{sec:num}, some numerical experiments are carried out to demonstrate the proposed theories. Finally, some conclusions can be found in \cref{con}.

\section{Preliminaries}\label{sec:pre}                              
\begin{figure}[!ht]
\centering
\begin{tikzpicture}[scale=1.5]
\draw[step=0.25, gray, thin] (-0.4, -0.4) grid (4.4, 4.4);
\draw[step=1.0, black, very thick] (-0.4, -0.4) grid (4.4, 4.4);
\foreach \x in {0,...,4}
\foreach \y in {0,...,4}{
\fill (1.0 * \x, 1.0 * \y) circle (1.5pt);
}
\fill[brown, opacity=0.4] (1.0, 1.0) rectangle (2.0, 2.0);
\node at (1.5, 1.5) {$K_j$};
\draw [dashed, very thick, fill=gray, opacity=0.6] (0, 0) rectangle (3, 3);
\node[above right] at (2.25, 2.25) {$K_j^1$};
\draw [dashed, very thick, fill=cyan, opacity=0.5] (3.25, 1.25) rectangle (3.5, 1.5);
\node at (3.375, 1.375) {$h$};
\end{tikzpicture}
\caption{An illustration of the two-scale mesh, a fine element $h$, a coarse element $K_j$ and its oversampling coarse element $K_j^m$ with the oversmapling layer $m=1$.}
\label{fig:grid}
\end{figure}
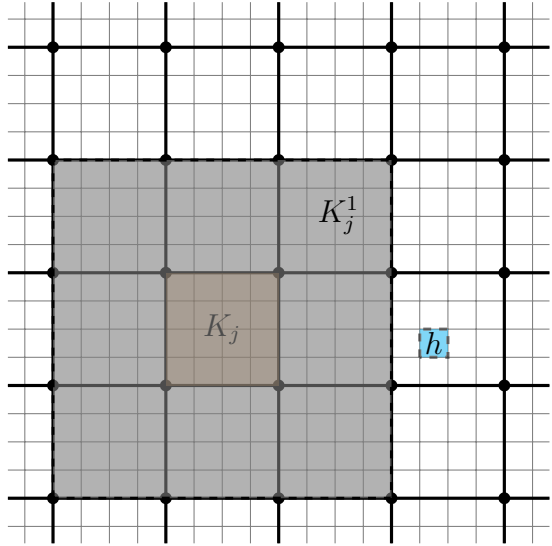

We consider the following Helmholtz equation for heterogeneous media in the bounded space domain $\Omega\subset\Real^d$ where $d=2$ or $3$:
\begin{equation}\label{eq:ell1}
\left\{
\begin{aligned}
-\nabla\cdot(A\nabla u)-k^2u=&f \quad &&\text{in} \quad \Omega, \\
u=&0\quad &&\text{on}\quad\Gamma_D,\\
A\nabla u\cdot \bm{n}-\mathrm{i}ku =&0 \quad &&\text{on}\quad \Gamma_R,
\end{aligned}
\right.
\end{equation}
where $\Gamma=\Gamma_\mathup{D}\cup\Gamma_\mathup{N}\cup\Gamma_\mathup{R}$ is a Lipschitz continuous boundary where  $\Gamma_\mathup{D}, \Gamma_\mathup{N},\Gamma_\mathup{R}$ represents the Dirichlet, Neumann and Robin boundary conditions respectively,  $\bm{n}$ is the unit outward normal vector to the boundary, $k\in\Real$ is a positive wavenumber, $ f
 \in L^2(\Omega)$ represents a harmonic source,  $\mathrm{i}$ denotes the imaginary unit and the scalar diffusion coefficient $A$ is a piecewise constant with respect to a quadrilateral background mesh $\mathcal{T}_{\varepsilon}$ with mesh size $O(\varepsilon)$ and $0 <\varepsilon< 1$. On each quadrilateral, $A$ takes either the value $\varepsilon^2$ or 1.
We define the Sobolev space $V=\CurlyBrackets*{u\in H^1(\Omega) : u=0 \, \text{on} \,\Gamma_\mathup{D}}$ and the $k$-weighted norm 
\[
\Norm{u}_V:=\Norm{u}_{1,A,k}:=\sqrt{k^2\|u\|^2+\|A^{1/2}\nabla u\|^2},
\]
where $\Norm{\cdot}$ denotes the $L^2$-norm over $\Omega$. We write the boundary value problem \cref{eq:ell1} in a variational form and find a solution $u\in V$ such that for all $v\in V$,
\begin{equation}\label{eq:weak}
\begin{aligned}
\int_\Omega A\nabla u\cdot\nabla\overline{v} \dx x-\mathrm{i}k\int_{\Gamma}u \overline{v}\dx s-\int_{\Omega}k^2 u\overline{v}\dx x &=\int_{\Omega} f\overline{v} \dx x, &&\forall v\in {V},\\
\end{aligned}
\end{equation}
where $\dx s$ represents the element of arc length along boundary $\Gamma$ and $\overline{\cdot}$ is the complex conjugation. To simplify the notations, 
 the sesquilinear form $\mathcal{B}: V\times V\rightarrow\Complex$ satisfies
 \begin{equation}\label{bform}
\mathcal{B}(u,v) :=(A\nabla u, \nabla v)-\mathrm{i}k(u,v)_{\Gamma}-k^2(u,v),
 \end{equation}
where $\RoundBrackets*{u,v}=\int_{\Omega}u\cdot\overline{v}\dx x$ and $\RoundBrackets*{u,v}_{\Gamma}=\int_{\Gamma}u\cdot\overline{v}\dx s$. Then we rewrite \cref{eq:weak} in the following
\begin{equation}\label{vf}
\mathcal{B}(u,v)=(f,v).
\end{equation}
 The well-posedness of the weak form \cref{vf} can be found in \cite{Graham2018} and for all $u, v \in V$ that satisfy the following inf-sup condition is as follows
 \[
\inf_{u\in V}\sup_{v\in V}\frac{\Re\mathcal{B}(u,v)}{\Norm{u}_{1,A,k}\Norm{v}_{1,A,k}}\geq\gamma>0.
\]
Let $\mathcal{T}_H$ be a standard quadrilteralization of the domain $\Omega$ with the mesh size $H$, where we call $\mathcal{T}_H$ coarse grid, $H>0$ being the coarse grid size. We refer to this partition as coarse grids, and the produced elements as the coarse elements. In each coarse element $K\in\mathcal{T}_H$, $K$ is further partitioned into a union of connected fine grid blocks. We denote the fine-grid partition as $\mathcal{T}_h$ with $h>0$ being the fine grid size. For each $K_j\in \mathcal{T}_H$ with $1\leq j\leq N$, $N$ is the number of the coarse grid elements. As is shown in \cref{fig:grid} we define an oversampled domain $K_j^m(m\geq 1)$ in the following 
$$
K_j^m: =\text{int}\CurlyBrackets*{\bigcup_{K\in\mathcal{T}_H}\CurlyBrackets*{  \overline{K_j^{m-1}}\bigcap\overline{K}\neq \varnothing}\bigcup\overline{K_j^{m-1}}},
$$
where $\text{int}(S)$ and $\overline{S}$ represent the interior and the closure of a set $S$, and the initial value $K_j^0=K_j$ for each element.
\section{The construction of the CEM-GMsFEM basis function}\label{sec:cemgms}
In this section, We now present our CEM-GMsFEM to efficiently solve \cref{eq:ell1}. The main process is as the following: we first construct auxiliary basis functions by solving an eigenvalue problem in an coarse element $K_j$ , then move on the auxiliary basis functions to the multiscale basis functions throughout oversampling areas by utilizing the idea of constrained energy minimization , finally the computational method will also been given and we solve Helmholtz problems in the newly obtained multiscale basis functions spaces.

Let $V(K_j)$ be the snapshot space on each coarse grid block $K_j$, and we use the method of the spectral problem to solve an eigenvalue problem on $K_j$: find eigenvalues $\lambda_j^{i}\in \mathbb{R}$ and basis functions $\phi_j^{i}\in H^1(K_j)$ such that for all $v\in V(K_j)$,
\begin{equation}\label{eig}
a_j(\phi_j^i,v)=\lambda_j^i s_j(\phi_j^i,v), \quad\forall v\in {H^1(K_j)},
\end{equation}
where 
$$
\begin{aligned}
a_j(\phi_j^i,v)=&\int_{K_j}A\nabla \phi_j^i\cdot\nabla\overline v \dx x,\\
s_j(\phi_j^i,v)=&\int_{K_j}{\tilde{A}(x)}\phi_j^i\overline{v}\dx x,\quad \tilde{A}(x) :=\sum_{k=1}^{N_v}A\nabla\eta_{j,k}^{i}\cdot\nabla\eta_{j,k}^{i},
\end{aligned}
$$
where, $N_v$ is the number of vertices contained in an element, to be specific, $N_v=4$ for a quadrilateral mesh and  $\CurlyBrackets*{\eta_{j,1}^{i},\eta_{j,2}^{i},\cdots,\eta_{j,N_v}^{i}}$ is the set of Lagrange basis functions on the coarse element $K_j\in\mathcal{T}_H.$  For the function $\tilde{A}(x)$, there is a non-negative constant $c_*$ such that 
$$c_*^{-1}H^{-2}\min(1,\varepsilon^2)\leq \tilde{A}(x)\leq c_* H^{-2}\max(1,\varepsilon^2).$$
Let the eigenvalues $\lambda_j^i$ in the ascending order:
$$
0=\lambda_j^0<\lambda_j^1\leq\lambda_j^2\leq\cdots\lambda_j^{l_j+1}\leq\cdots ,
$$ 
and we use the first $l_j$ eigenvalue functions corresponding to the eigenvalues to construct the local auxiliary space  $V_{\mathup{aux}}^{j}=\CurlyBrackets{\phi_j^1,\phi_j^2,\cdots,\phi_j^{l_j}}.$ The global auxiliary space $V_{\mathup{aux}}$ is the sum of these local auxiliary spaces, namely 
$V_{\mathup{aux}}= \bigoplus_{j=1}^N V_{\mathup{aux}}^{j}$, which will be used to construct multiscale basis functions.
The next we give the definition of the so called $\phi_j^i$-orthogonal, for a given a function $\phi_j^i\in V_{\mathup{aux}}$, $\psi\in V$, and we define
$$
s(\phi_j^i,\psi)=1, \quad s(\phi_{j'}^{i'},\psi)=0 \,\,\text{if}\, j'\neq j\,\text{or}\,i'\neq i.
$$
Based on the $\phi_j^i$-orthogonal, we can obtain that for any $v\in V$
$$s(\phi_j^i,v)=\sum_{j=1}^{N}s_j(\phi_j^i,v).$$
The orthogonal projection $\pi_j$ from $V(K_j)$ onto $V_\mathup{aux}^{j}$ is
$$\pi_{j}(v) :=\sum_{i=1}^{l_j}\frac{s(\phi_j^i,v)}{s(\phi_j^i,\phi_j^i)}\phi_j^i, \quad \forall v\in V(K_j),$$
and the global projection is $\pi :=\sum_{j=1}^{N}\pi_j$ from $L^2$ to $V_\mathup{aux}$
\footnote{We use a zero-extension here, which extends each $V_{\mathup{aux}}^j$ into $L^2(\Omega)$. }
We can immediately derive the following \cref{inter}, which shows an important property of the global projection $\pi$.
\begin{lemma}\label{inter}
In each $K_j\in\mathcal{T}_H$, for all $v\in H^1(K_j)$,
\begin{equation}\label{pi1}
\Norm{v-\pi_j v}_{s(K_j)}^2\leq \frac{\Norm{v}_{a(K_j)}^2}{\lambda_{l_j+1}^j}\leq\Lambda^{-1}\Norm{v}_{a(K_j)}^2.
\end{equation}
where $\Lambda=\max_{1\leq j\leq N}\lambda_j^{l_j+1}$, and  
\begin{equation}\label{pi2}
\norm{\pi_jv}_{s(K_j)}^2=\Norm{v}_{s(K_j)}^2-\Norm{v-\pi_jv}_{s(K_j)}^2\leq	\Norm{v}_{s(K_j)}^2.
\end{equation}
\end{lemma}
\subsection{Multiscale basis functions}
In order to deal with the lack of hermitivity of the $\mathcal{B}$, we need to define two bounded operators $T_m=\sum_{K_j\in\mathcal{T}_H}T_{j,m}$ and $T_m^*=\sum_{K_j\in\mathcal{T}_H}T_{j,m}^*$ both from $L^2(\Omega) $ to $H^1(\Omega)$ to construct the test space and trial space of the following variational problems.  For each coarse element $K_j\in\mathcal{T}_H$ and its oversampling domain $K_j^m\subset\Omega$ by enlarging $K_j$ for $m$ coarse grid layers, we define the multiscale basis function $T_{j,m}\psi_{j,m}^i\in V_0(K_j^m)$, find $T_{j,m}\psi_{j,m}^i\in V_0(K_j^m)$ such that
\begin{equation}\label{msbasis1}
\mathcal{B}(T_{j,m}\psi_{j,m}^i,v)+s(\pi T_{j,m}\psi_{j,m}^i,\pi v)=s(\pi_j\psi_{j,m}^i,\pi v),\quad\forall v\in V_0(K_j^m),
\end{equation}
where $V_0(K_j^m)$ is the subspace of $V(K_j^m)$ with zero trace on $\partial K_j^m$ and $V(K_j^m)$ is the restriction of $V$ in $K_j^m.$ 
Now, our multiscale finite element space $V_{\mathup{ms}}$ can be defined by solving an variational problem \cref{msbasis1}
$$
V_{\mathup{ms}}=\mathup{span}\CurlyBrackets*{T_{j,m}\psi_{j,m}^{i}\,|\, 1\leq i\leq l_j,\,1\leq j\leq N}.
$$
The global multiscale basis function $T_j\psi_j^i\in V$ is defined in a similar way, 
\begin{equation}\label{msbasis2}
\mathcal{B}(T_j\psi_{j}^i,v)+s(\pi T_j\psi_{j}^i,\pi v)=s(\pi_j\psi_{j}^i,\pi v),\quad\forall v\in V(K_j),
\end{equation}
and $T=\sum_{K_j\in\mathcal{T}_H}T_{j}$. Thereby, the global multiscale finite element space $V_{\mathup{glo}}$  is defined by
$$
V_{\mathup{glo}}=\mathup{span}\CurlyBrackets*{T_j\psi_j^i\,|\, 1\leq i\leq l_j,\,1\leq j\leq N}.
$$
Similarly, for the local operator $T_{j,m}^*$ from $L^2(K_j^m)$ to $H^1(K_j^m)$,
\begin{equation}\label{msbasis3}
\mathcal{B}(v,T_{j,m}^*\psi_{j,m}^i)+s(\pi v,\pi T_{j,m}^*\psi_{j,m}^i)=s(\pi v,\pi_j\psi_{j,m}^i),\quad\forall v\in V_0(K_j^m),
\end{equation}
where $T_{j,m}^*(v)=\overline{T_{j,m}(\overline{v})}$.
\begin{remark}
Due to the operator $\mathcal{B}$ is in the complex domain, there is no minimum solutions of the minimization problem listed in \cite{Chung2018}.
\end{remark}
 Now, another multiscale finite element space $V_{\mathup{ms}}^*$ can be defined by solving the above variational problem \ref{msbasis3}
$$
V_{\mathup{ms}}^*=\mathup{span}\CurlyBrackets*{T_{j,m}^*\psi_{j,m}^{i}\,|\, 1\leq i\leq l_j,\,1\leq j\leq N}.
$$
The global multiscale basis function $T_j^*\psi_j^i\in V$ is defined in a similar way,
\begin{equation}\label{msbasis4}
\mathcal{B}(v,T_j^*\psi_{j}^i)+s(\pi v,\pi T_j^*\psi_{j}^i)=s(\pi v,\pi_j\psi_{j}^i),\quad\forall v\in V(K_j),
\end{equation}
where $T^*=\sum_{K_j\in\mathcal{T}_H}T_{j}^*$. Thereby, another global multiscale finite element space $V_{\mathup{glo}}^*$  is defined by
$$
V_{\mathup{glo}}^*=\mathup{span}\CurlyBrackets*{T_j^*\psi_j^i\,|\, 1\leq i\leq l_j,\,1\leq j\leq N}.
$$
The existence of the above two variational problems will be shown in the analysis part. In the following, we use $V_{\mathup{ms}}$ and  $V_{\mathup{ms}}^*$ as the new test space and trial space of Petrov-Galerkin method to find the approximated solution of  \cref{vf}: Find $u_{\mathup{ms}}\in V_{\mathup{ms}}$ such that 
\begin{equation}\label{cemvar}
\mathcal{B}(u_{\mathup{ms}},v)=(f,v),\quad \forall v\in V_{\mathup{ms}}^*.
\end{equation}
For further analysis, we need to clarify an important orthogonality property of the global multiscale finite element space $V_{\mathup{glo}}$ and $V_{\mathup{glo}}^*$.
\begin{lemma}\label{lm3.2}
The space $V$ is decomposed as
$$V=\widetilde{V}\oplus V_{\mathup{glo}}^*,$$
where $\widetilde{V}=\CurlyBrackets*{v\in V\,|\, \pi(v)=0}$ and $W=\CurlyBrackets*{v\in V\,|\, \mathcal{B}(v,w)=0,\forall w\in V_{\mathup{glo}}^*}$, then $\widetilde{V}=W$.
\end{lemma}
\begin{proof} By using the fact in \cref{msbasis4}, it is easy to see $\tilde V\subset W$.  For another direction, $\forall v_0\in W$, we have $\mathcal{B}(v_0,w)=0$. Then we can obtain $s(\pi v_0,\pi w)=s(\pi v_0,\pi\psi_j^i)$. Due to the arbitrary $w\in V_{\mathup{glo}}^*$, we can get that $\pi v_0=0$ and $W\subset \tilde V$.
\end{proof} 
\begin{lemma}\label{lm3.3}
The space $V$ is also decomposed as
$$V=\widetilde{V}\oplus V_{\mathup{glo}},$$
where $\widetilde{V}=\CurlyBrackets*{v\in V\,|\, \pi(v)=0}$ and $\widetilde{W}=\CurlyBrackets*{v\in V\,|\, \mathcal{B}(w,v)=0,\forall w\in V_{\mathup{glo}}}$, then $\widetilde{V}=\widetilde{W}$.
\end{lemma}
\begin{proof} By using the fact in \cref{msbasis2}, it is easy to see $\tilde V\subset\widetilde{W}$.  For another direction, $\forall v_0\in \widetilde{W}$, we have $\mathcal{B}(w,v_0)=0$. Then we can obtain $s(\pi w,\pi v_0)=s(\pi v_0,\pi\psi_j^i)$. Due to the arbitrary $w\in V_{\mathup{glo}}$, we can get that $\pi v_0=0$ and $\widetilde{W}\subset \tilde V$.
\end{proof}
\begin{remark}
The above two lemmas show a relationship of “orthogonality” between $\widetilde{V}$ and its complementary space $V_{\mathup{glo}}$,  $V_{\mathup{glo}}^*$ concerning the bilinear form $\mathcal{B}(\cdot, \cdot).$ However, we must be cautious in using the term “orthogonal” since $\mathcal{B}(\cdot, \cdot)$ cannot define an inner product on $V$.
\end{remark}
\section{Analysis}\label{sec:anal}
In this section, we give the stability and convergence of CEM-GMsFEM by using the global multiscale basis functions. We firstly prove that the sesquilinear form of $a$ is coercive, continuous and bounded.  Then we show the exponential decay property of the multiscale basis functions. In particular, we will show that the the global basis function and the corresponding local basis function  are the same the if the oversampling region is sufficiently large. Finally, we prove the convergence of the multiscale solution in our main result  \cref{thm:main}. Here,  the approximated solution $u_{\mathup{glo}}\in V_{\mathup{glo}}$ obtained in the global multiscale space $V_{\mathup{glo}}$ is defined by
\begin{equation}\label{glo}
\mathcal{B}(u_{\mathup{glo}},v)=(f,v),\quad\forall v\in V_{\mathup{glo}}^*.
\end{equation}
\begin{lemma}\label{lm1}
If the mesh size $H$, the wave number $k$, and diffusion parameter $\varepsilon$ satisfies the resolution condition such that
$$\sqrt{2c_*}k H\varepsilon^{-1}\leq 1,$$
the sesquilinear form defined in \cref{bform} is  continuous on $V$ , that is
\begin{equation}\label{lem2.1}
\abs{\mathcal{B}(u,v)}\leq c\Norm{u}_{1,A,k}\Norm{v}_{1,A,k}\quad \forall\, u, v\in V,
\end{equation}
where $c$ is independent of $k$, and the sesquilinear form is coercive on $ \tilde{V}$ defined in \cref{lm3.2},
\begin{equation}\label{lem2.2}
\epsilon_0\Norm{v}_{a}^2\leq\Re\mathcal{B}(v,v),\quad \forall v\in \tilde{V},
\end{equation}
where $\mathup{Re}$ and $\mathup{Im}$ represents the real and imaginary parts of a complex number, respectively.
\end{lemma}
\begin{proof}
Rewrite the definition of $\mathcal{B}(u,v)$ in \cref{bform}, we have 
$$
\mathcal{B}(u,v)=b_1+b_2+b_{\Gamma},
$$
where 
\begin{equation*}
\begin{aligned}
b_1(u,v) :&= \int_{\Omega} \RoundBrackets*{A\nabla u\cdot\nabla\overline{v}+k^2u\overline{v}}\dx x,\\
\quad b_2(u,v) :&= -2k^2\int_{\Omega}u\overline{v}\dx x,
\quad b_{\Gamma} :=-\mathrm{i}k\int_{\Gamma}u\overline{v}\dx x.    
\end{aligned}
\end{equation*}
By using the Cauchy--Schwartz inequality we obtain 
$$
\begin{aligned}
\abs{b_1(u,v)}\leq &\Norm{u}_{1,A,k}\Norm{v}_{1,A,k},\\
\abs{b_2(u,v)}\leq &2\Norm{u}_{1j,A,k}\Norm{v}_{1,A,k},\\
\end{aligned}
$$
and the trace inequalities from \cite{Grisvard2011},
 $$  
\abs{b_{\Gamma}(u,v)}\leq C_{\Gamma}\Norm{u}_{1,A,k}\Norm{v}_{1,A,k},\\
$$
where $C_{\Gamma}$ is independent of  $k$. 
By using the interpolation properties of the operator $\pi$ in \cref{inter}  and assumption about the resolution condition
$$
\Re\mathcal{B}(v,v)=\norm{v}_{a}^2-k^2\Norm{v}^2
\geq (\Lambda-c_*k^2H^2\varepsilon^{-2})\Norm{v}_{s}^2\geq\epsilon_0\norm{v}_{a}^2,
$$
where the eigenvalue $\Lambda$ is chosen such that $\Lambda\geq\frac{1}{2}+\epsilon_0$ and $\epsilon_0$ is a small positive value.
\end{proof}
The following lemma gives the well-posedness of the global multiscale problem \cref{glo} by proving the inf-sup condition.
\begin{lemma}\label{thm:well}
The bilinear form $\mathcal{B}$ satisfies the following inf-sup condition: there exists a constant $\gamma>0$ depends on $k$ such that
\[
\inf_{u_{\mathup{glo}}\in V_{\mathup{glo}}}\sup_{u_{\mathup{glo}^*\in V_{\mathup{glo}}^*}}\frac{\Re\mathcal{B}(u_{\mathup{glo}},u_{\mathup{glo}}^*)}{\Norm{u_{\mathup{glo}}}_{1,A,k}\Norm{u_{\mathup{glo}}^*}_{1,A,k}}\geq\gamma>0.
\]
\end{lemma}
Based on the previous decomposition of the space $V$, the Fortin trick \cite{Boffi2013} suggests that we only need to check 
\[
\inf_{v_1\in \widetilde{V}}\sup_{v_2\in\widetilde{V}}\frac{\Re\mathcal{B}(v_1,v_2)}{\Norm{v_1}_{1,A,k}\Norm{v_2}_{1,A,k}}\geq\gamma>0.
\]
\begin{proof}
For any $v_1\in\widetilde{V}$, we apply the inf-sup condition of the $\mathcal{B}(u,v)$ in $V\times V$, for any $v_2\in V$, there exists a $\phi\in V$ and $\Norm{\phi}_{1,A,k}=1.$ such that $$\mathup{Re}\mathcal{B}(v_1,\phi)\geq\gamma\Norm{v_1}_{1,A,k}.$$ 
Then we can define $v_2=\phi-T^*\phi\in\widetilde{V}$, and
$\mathcal{B}(v_1,v_2)=\mathcal{B}(v_1,\phi)-\mathcal{B}(v_1,T^*\phi)=\mathcal{B}(v_1,\phi),$
where the last equality comes from the fact that $\mathcal{B}(v_1,T^*\phi)=0$. Since $\mathcal{B}(v_1,T^*\phi)=\mathcal{B}(T\overline{\phi},\overline{v_1})=0$ where $T\overline{\phi}\in V_{\mathup{glo}}$ and $\overline{v_1}\in \widetilde{V}$, by using of the orthogonality in \cref{lm3.3}, especially, we can obtain $\mathcal{B}(T\overline{\phi},\overline{v_1})=0$. 
Finally we can obtain the desired inequality that
$$\mathup{Re}\mathcal{B}(v_1,v_2)=\mathup{Re}\mathcal{B}(v_1,\phi)\geq\gamma\Norm{v_1}_{1,A,k}\geq\gamma\Norm{v_1}_{1,A,k}\Norm{v_2}_{1,A,k},
$$
where the last inequality comes from that $\Norm{v_2}_{1,A,k}\leq \Norm{\phi}_{1,A,k}+\Norm{T^*\phi}_{1,A,k}\leq (1+C)\Norm{\phi}_{1,A,k}$ due to the bounded operator $T^*$.
\end{proof}
After obtaining the well-posedness of the global problem, an error estimate of the global solution can be derived in \cref{lm2} .
\begin{lemma}\label{lm2}
By using the resolution condition in \cref{lm1}, let $u_{\mathup{glo}}$ be the solution of \cref{glo} and $u$ be the real solution of the problem \cref{vf}. Then 
$$
\Norm{u_{\mathup{glo}}-u}_{a}\leq\frac{1}{\epsilon_0\sqrt{\Lambda}}\Norm{f}_{s^{-1}},
$$
where $\Lambda=\max_{1\leq j\leq N}\lambda_{l_j+1}^j$, and
$$\Norm{f}_{s^{-1}}:=\sup_{v\in L^2(\Omega), v\neq 0}\frac{\int_{\Omega}fv\dx x}{\Norm{v}_s}.$$
\end{lemma}
\begin{proof}
By direct computation in \cref{glo} and \cref{vf}, we obtain the Galerkin orthogonality that for any $v\in V_{\mathup{glo}}^*$ 
$$\mathcal{B}\RoundBrackets*{u-u_{\mathup{glo}},v}=0.$$
In \cref{lm3.2}, we have $\mathcal{B}\RoundBrackets*{u-u_{\mathup{glo}},v}=0 $, thus $u-u_{\mathup{glo}}\in \tilde{V}$ and $\pi(u-u_{\mathup{glo}})=0$. By using coercivity of $\mathcal{B}$ on $\tilde{V}$ and variational form in \cref{vf}, we have 
\begin{align*}
\epsilon_0\Norm{u-u_{\mathup{glo}}}_{a}^2\leq&\Re\mathcal{B}(u-u_{\mathup{glo}},u-u_{\mathup{glo}})\\
\leq&\abs{\mathcal{B}(u,u-u_{\mathup{glo}})}\\
=&\abs{(f,u-u_{\mathup{glo}})}\\
\leq&\Norm{f}_{s^{-1}}\Norm{u-u_{\mathup{glo}}}_s=\Norm{f}_{s^{-1}}\Norm{u-u_{\mathup{glo}}-\pi(u-u_{\mathup{glo}})}_s\\
\leq&\frac{1}{\sqrt{\Lambda}}\Norm{f}_{s^{-1}}\Norm{u-u_{\mathup{glo}}}_a,
\end{align*}
where the last inequality comes from \cref{eig} directly.
\end{proof}
The following \cref{lm3} will  show the multiscale basis functions have the exponentially decaying property. Before we give the detailed proof, we need to define the cutoff functions with respect to these oversampling domains in the following. 

For the cutoff functions, in each $K_j$, let  $V_H$ be  the Lagrange basis function space of  $\mathcal{T}_H$, we define $\eta_j^{n,m}\in V_H, m>n$ such that $0\leq\eta_j^{n,m}\leq 1$ and 
$$
\begin{aligned}
\eta_j^{n,m}=0,&\quad\text{in}\,K_{j,n},\\
\eta_j^{n,m}=1,&\quad\text{in}\,\Omega\setminus K_{j,m},\\
0\leq\eta_j^{n,m}\leq1,&\quad\text{in}\,K_{j,m}\setminus K_{j,n}.\\
\end{aligned}
$$
\begin{lemma}\label{lm3}
There exists a constant $0<\beta<1$, independent of $H, m$ and $k$ such that for any $t_j\in V_H$,
$$
\Norm{T_jt_j}_{a(\Omega\setminus{K_{j,m}})}^2+\Norm{\pi T_jt_j}_{s(\Omega\setminus K_{j,m})}^2\leq\beta^m\RoundBrackets*{\Norm{T_jt_j}_{a}^2+\Norm{\pi T_jt_j}_{s}^2} ,
$$
where $\beta=\frac{C}{1+C}.$
\end{lemma}
\begin{proof}
Let $t_j\in V_H$, and choose $w=\eta_j^{m-1,m}T_jt_j$ with support only outside $K_{j,m-1}$, therefore we obtain $\mathcal{B}(T_jt_j,w)+s(\pi T_jt_j,\pi w)=0$. By recalling the definition of the operator $\mathcal{B}$ we have 
$$
\begin{aligned}
&\int_{\Omega}A\nabla T_jt_j\cdot\nabla\overline{\RoundBrackets*{\eta_j^{m-1,m}T_jt_j}}\dx  x +\int_{\Omega}\tilde{A}\pi T_jt_j\cdot\pi\overline{\RoundBrackets*{\eta_j^{m-1,m}T_jt_j}}\dx x\\
=&\mathrm{i}k\int_{\Gamma} T_jt_j\cdot\overline{\RoundBrackets*{\eta_j^{m-1,m}T_jt_j}}\dx s +k^2\int_{\Omega}T_jt_j\cdot\overline{\RoundBrackets*{\eta_j^{m-1,m}T_jt_j}}\dx x.\\
\end{aligned}
$$
By using the properties of the cutoff functions to observe that $1-\eta_j^{m-1,m}=0$ in $\Omega\setminus K_{j,m}$ and $1-\eta_j^{m-1,m}=1$ in $K_{j,m-1}$ , we can obtain
$$
\begin{aligned}
&\Norm{T_jt_j}_{a(\Omega\setminus{K_{j,m}})}^2+\Norm{\pi T_jt_j}_{s(\Omega\setminus K_{j,m})}^2\\
=&-\int_{K_{j,m}\setminus K_{j,m-1}}\eta_j^{m-1,m}A\nabla T_jt_j\cdot\nabla\overline{T_jt_j}\dx x-\int_{K_{j,m}\setminus K_{j,m-1}}\overline{ T_jt_j}A\nabla T_jt_j\cdot\nabla\eta_j^{m-1,m}\dx x\\
&-\int_{K_{j,m}\setminus K_{j,m-1}}\tilde{A}\pi T_jt_j\cdot\pi\overline{\RoundBrackets*{\eta_j^{m-1,m}T_jt_j}}\dx x\\
&+\mathrm{i}k\int_{\Gamma} T_jt_j\cdot\overline{\RoundBrackets*{\eta_j^{m-1,m}T_jt_j}}\dx s+k^2\int_{\Omega}T_jt_j\cdot\overline{\RoundBrackets*{\eta_j^{m-1,m}T_jt_j}}\dx x\\
=&-\Re\int_{K_{j,m}\setminus K_{j,m-1}}\eta_j^{m-1,m}A\nabla T_jt_j\cdot\nabla\overline{T_jt_j}\dx x-\Re\int_{K_{j,m}\setminus K_{j,m-1}} \overline{T_jt_j}A\nabla T_jt_j\cdot\nabla\eta_j^{m-1,m}\dx x\\
&-\Re\int_{K_{j,m}\setminus K_{j,m-1}}\tilde{A}\pi T_jt_j\cdot\pi\overline{\RoundBrackets*{\eta_j^{m-1,m}T_jt_j}}\dx x\\
&+\Re\RoundBrackets*{\mathrm{i}k\int_{\Gamma} T_jt_j\cdot\overline{\RoundBrackets*{\eta_j^{m-1,m}T_jt_j}}\dx s}+\Re\RoundBrackets*{k^2\int_{\Omega}T_jt_j\cdot\overline{\RoundBrackets*{\eta_j^{m-1,m}T_jt_j}}\dx x}\\
\leq&\abs{\Re\int_{K_{j,m}\setminus K_{j,m-1}}\eta_j^{m-1,m}A\nabla T_jt_j\cdot\nabla\overline{T_jt_j}\dx x}+\abs{\Re\int_{K_{j,m}\setminus K_{j,m-1}} \overline{T_jt_j}A\nabla T_jt_j\cdot\nabla\eta_j^{m-1,m}\dx x}\\
&+\abs{\Re\int_{K_{j,m}\setminus K_{j,m-1}}\tilde{A}\pi T_jt_j\cdot\pi\overline{\RoundBrackets*{\eta_j^{m-1,m}T_jt_j}}\dx x}\\
&+\abs{\Re\RoundBrackets*{\mathrm{i}k\int_{\Gamma} T_jt_j\cdot\overline{\RoundBrackets*{\eta_j^{m-1,m}T_jt_j}}\dx s}}+\abs{\Re\RoundBrackets*{k^2\int_{\Omega}T_jt_j\cdot\overline{\RoundBrackets*{\eta_j^{m-1,m}T_jt_j}}\dx x}}\\
=&\sum_{i=1}^5I_i.
\end{aligned}
$$
For $I_4$, due to the property of the cutoff function, we can obtain
\begin{equation*}
\begin{aligned}
I_4=\abs{\text{Re}\RoundBrackets*{\mathrm{i}k\int_{\Gamma} T_jt_j\cdot\overline{\RoundBrackets*{\eta_j^{m-1,m}T_jt_j}}\dx s}}
=0.
\end{aligned}
\end{equation*}
Combining with the properties of $\eta_j^{m-1,m}$ and the resolution condition gives
\begin{equation*}\label{eq:k}
\begin{aligned}	
I_5\leq\abs{k^2\int_{\Omega\setminus K_{j,m-1}}T_jt_j\cdot\overline{T_jt_j}\dx x}
\leq \frac{1}{2}\RoundBrackets{\Norm{T_jt_j}_{s(\Omega\setminus K_{j,m})}^2+\Norm{T_jt_j}_{s(K_{j,m}\setminus K_{j,m-1})}^2}.
\end{aligned}
\end{equation*}
Also, \cref{pi1} provides an estimate for $\Norm{T_jt_j}_{s(K_{j,m}\setminus K_{j,m-1})}$ and $\Norm{T_jt_j}_{(\Omega\setminus K_{j,m})}$, 
$$
\begin{aligned}
\Norm{T_jt_j}_{s(K_{j,m}\setminus K_{j,m-1})}^2&\leq \Norm{T_jt_j-\pi T_jt_j}_{s(K_{j,m}\setminus K_{j,m-1})}^2+\Norm{\pi T_jt_j}_{s(K_{j,m}\setminus K_{j,m-1})}^2\\
&\leq \frac{1}{\Lambda}\Norm{T_jt_j}_{a(K_{j,m}\setminus K_{j,m-1})}^2+\Norm{\pi T_jt_j}_{s(K_{j,m}\setminus K_{j,m-1})}^2,
\end{aligned}
$$
and
$$
\begin{aligned}
\Norm{T_jt_j}_{s(\Omega\setminus K_{j,m})}^2&\leq \Norm{T_jt_j-\pi T_jt_j}_{s(\Omega\setminus K_{j,m})}^2+\Norm{\pi T_jt_j}_{s(\Omega\setminus K_{j,m})}^2\\
&\leq \frac{1}{\Lambda}\Norm{T_jt_j}_{a(\Omega\setminus K_{j,m})}^2+\Norm{\pi T_jt_j}_{s(\Omega\setminus K_{j,m})}^2.
\end{aligned}
$$
By the definition of $\tilde{A}$, it is easy  to show
$$ A(x)\nabla\eta_j^{m-1,m}\cdot \nabla \eta_j^{m-1,m}\leq \tilde{A}(x),$$
and using the properties of $\eta_j^{m-1,m}$, the estimates of $I_2$ follows
$$I_2\leq\Norm{T_jt_j}_{a(K_{j,m}\setminus K_{j,m-1})}\Norm{T_jt_j}_{s(K_{j,m}\setminus K_{j,m-1})}.$$
For the term $I_3$, using the Cauchy-Schwarz inequality, we have 
$$
\begin{aligned}
I_3 \leq&\Norm{\pi T_jt_j}_{s(K_{j,m}\setminus K_{j,m-1})}\Norm{\pi\RoundBrackets{\eta_j^{m-1,m}T_jt_j}}_{s(K_{j,m}\setminus K_{j,m-1})}\\
\leq&\Norm{\pi T_jt_j}_{s(K_{j,m}\setminus K_{j,m-1})}\Norm{\eta_j^{m-1,m}T_jt_j}_{s(K_{j,m}\setminus K_{j,m-1})}\\
\leq&\Norm{\pi T_jt_j}_{s(K_{j,m}\setminus K_{j,m-1})}\Norm{T_jt_j}_{s(K_{j,m}\setminus K_{j,m-1})}.
\end{aligned}
$$
For the last term $I_1$, Due to the definition of  $\eta_j^{m-1,m}$, we can derive
$$I_1\leq\Norm{T_jt_j}_{a(K_{j,m}\setminus K_{j,m-1})}^2.$$
Combining all above estimates of $I_1$ to $I_5$ together,  we obtain
$$
\begin{aligned}
&\RoundBrackets*{1-\frac{1}{2\Lambda}}\Norm{T_jt_j}_{a(\Omega\setminus{K_{j,m}})}^2+\frac{1}{2}\Norm{\pi T_jt_j}_{s(\Omega\setminus K_{j,m})}^2\\
\leq& C(\Norm{T_jt_j}_{a( K_{j,m}\setminus K_{j,m-1})}^2+\Norm{\pi T_jt_j}_{s(K_{j,m}\setminus K_{j,m-1})}^2), 
\end{aligned}
$$
also, we need to choose $\Lambda>\frac{1}{2}$ to future derive 
$$
\begin{aligned}
&\Norm{T_jt_j}_{a(\Omega\setminus{K_{j,m}})}^2+\Norm{\pi T_jt_j}_{s(\Omega\setminus K_{j,m})}^2\\\leq& C(\Lambda)(\Norm{T_jt_j}_{a( K_{j,m}\setminus K_{j,m-1})}^2+\Norm{\pi T_jt_j}_{s(K_{j,m}\setminus K_{j,m-1})}^2).
\end{aligned}
$$
Because of 
$$
\Norm{T_jt_j}_{a(\Omega\setminus K_{j,m})}^2+\Norm{T_jt_j}_{a( K_{j,m}\setminus K_{j,m-1})}^2=\Norm{T_jt_j}_{a(K_{j,m-1})}^2,
$$
$$
\Norm{\pi T_jt_j}_{s(\Omega\setminus K_{j,m})}^2+\Norm{\pi T_jt_j}_{s( K_{j,m}\setminus K_{j,m-1})}^2=\Norm{\pi T_jt_j}_{s(K_{j,m-1})}^2,
$$
then 
$$
\begin{aligned}
&\Norm{T_jt_j}_{a(\Omega\setminus K_{j,m})}^2+\Norm{\pi T_jt_j}_{s(\Omega\setminus K_{j,m})}^2\\\leq&\frac{C(\Lambda)}{1+C(\Lambda)}\RoundBrackets{\Norm{T_jt_j}_{a( \Omega\setminus K_{j,m-1})}^2+\Norm{\pi T_jt_j}_{s(\Omega\setminus K_{j,m-1})}^2}. \end{aligned}
$$
Let $\beta=\frac{C(\Lambda)}{1+C(\Lambda)}<1$, and after repeat the above the inequality
$$
\Norm{T_jt_j}_{a(\Omega\setminus K_{j,m})}^2+\Norm{\pi T_jt_j}_{s(\Omega\setminus K_{j,m})}^2\leq \beta^m\RoundBrackets*{\Norm{T_jt_j}_{a}^2+\Norm{\pi T_jt_j}_{s}^2}. 
$$
\end{proof}
\begin{lemma}\label{lemma2} 	
Keep the notations same as \cref{lm3}, then for $t_j\in V_H$, there exists a constant $C(\Lambda)$ such that 
$$
\|(T_j-T_{j,m}) t_j\|_{a}^2+\|\pi(T_j-T_{j,m}) t_j\|_{s}^2\leq C(\Lambda)\beta^{m-1}(\|T_jt_j\|_{a}^2+\|\pi T_jt_j\|_{s}^2).
$$
\end{lemma}
\begin{proof} 
By using the previous definition of the operators $T_j$ and $T_{j,m}$, subtract each other we can obtain for $w=(1-\eta_j^{m-1,m})T_jt_j-T_{j,m}t_j\in V_{j,m}$ ,
\begin{equation*}\label{eq10}
\mathcal{B}\RoundBrackets*{T_jt_j-T_{j,m}t_j,w}+s\RoundBrackets*{\pi(T_jt_j-T_{j,m}t_j), \pi w}=0.
\end{equation*}
After reformulating the above equation,
\begin{equation*}
\begin{aligned}
&\mathcal{B}\RoundBrackets*{T_jt_j-T_{j,m}t_j,T_jt_j-T_{j,m}t_j}+s\RoundBrackets*{\pi(T_jt_j-T_{j,m}t_j),\pi(T_jt_j-T_{j,m}t_j)}\\
=&\mathcal{B}(T_jt_j-T_{j,m}t_j,\eta_j^{m-1,m}T_jt_j)+s(\pi(T_jt_j-T_{j,m}t_j),\pi(\eta_j^{m-1,m}T_jt_j)).
\end{aligned}
\end{equation*}
Recalling the definition of the operator $\mathcal{B}$, we obtain
\begin{equation*}\label{eq:18}
\begin{aligned}
&\Norm*{(T_j-T_{j,m})t_j}_{a}^2+\Norm*{\pi(T_j-T_{j,m})t_j}_{s}^2\\
=&a(T_jt_j-T_{j,m}t_j,\eta_j^{m-1,m}T_jt_j)+s(\pi(T_jt_j-T_{j,m}t_j),\pi(\eta_j^{m-1,m}T_jt_j))\\
&+\mathrm{i}k(T_jt_j-T_{j,m}t_j,w)_{\Gamma}+k^2(T_jt_j-T_{j,m}t_j, w)\\
=&\text{Re}\RoundBrackets*{a(T_jt_j-T_{j,m}t_j,\eta_j^{m-1,m}T_jt_j)}+\text{Re}\RoundBrackets*{s(\pi(T_jt_j-T_{j,m}t_j),\pi(\eta_j^{m-1,m}T_jt_j))}\\
&+\text{Re}\RoundBrackets*{\mathrm{i}k(T_jt_j-T_{j,m}t_j,w)_{\Gamma}}+\text{Re}\RoundBrackets*{k^2(T_jt_j-T_{j,m}t_j, w)}\\
\leq &\abs{\text{Re}\RoundBrackets*{a(T_jt_j-T_{j,m}t_j,\eta_j^{m-1,m}T_jt_j)}}+\abs{\text{Re}\RoundBrackets*{s(\pi(T_jt_j-T_{j,m}t_j),\pi(\eta_j^{m-1,m}T_jt_j))}}\\
&+\abs{\text{Re}\RoundBrackets*{\mathrm{i}k(T_jt_j-T_{j,m}t_j,w)_{\Gamma}}}+\abs{\text{Re}\RoundBrackets*{k^2(T_jt_j-T_{j,m}t_j, w)}}\\
=&\sum_{i=1}^4I_i.
\end{aligned}
\end{equation*}
For $I_3$, we use the trace inequality in \cite{brenner2008} to obtain
\begin{equation*}
\begin{aligned}
I_3\leq&\abs{\text{Re}\RoundBrackets*{\mathrm{i}k(T_jt_j-T_{j,m}t_j,w)_{\Gamma}}}\\
\leq&\abs{k(T_jt_j-T_{j,m}t_j,\eta_j^{m-1,m}T_jt_j)_{\Gamma}}\\
\leq&\Norm*{(T_j-T_{j,m})t_j}_{1,A,k}\Norm{\eta_j^{m-1,m}T_jt_j}_{1,A,k}\\
\leq&\RoundBrackets{\Norm*{(T_j-T_{j,m})t_j}_a^2+k^2\Norm*{(T_j-T_{j,m})t_j}^2}^{1/2}\RoundBrackets{\Norm{\eta_j^{m-1,m}T_jt_j}_a^2+k^2\Norm{\eta_j^{m-1,m}T_jt_j}^2}^{1/2}.
\end{aligned}
\end{equation*}
By using resolution condition
\begin{equation*}
\begin{aligned}
k^2\Norm{(T_j-T_{j,m})t_j}^2\leq c_*k^2H^2\varepsilon^{-2}\Norm*{(T_j-T_{j,m})t_j}_s^2\leq \frac{1}{2}\Norm*{(T_j-T_{j,m})t_j}_s^2,\\
k^2\Norm{\eta_j^{m-1,m}T_jt_j}^2\leq c_*k^2H^2\varepsilon^{-2}\Norm*{T_jt_j}_{s(\Omega\setminus K_{j,m-1})}^2\leq \frac{1}{2}\Norm*{T_jt_j}_{s(\Omega\setminus K_{j,m-1})}^2.
\end{aligned}
\end{equation*}
Furthermore, 
$$
\begin{aligned}
\Norm{(T_j-T_{j,m})t_j}_s^2&\leq\Norm{T_jt_j-T_{j,m}t_j-\pi(T_jt_j-T_{j,m}t_j)}_s^2+\Norm{\pi ((T_j-T_{j,m})t_j)}_s^2\\
&\leq\frac{1}{\Lambda}\Norm{(T_j-T_{j,m})t_j}_a^2+ \Norm{\pi((T_j-T_{j,m})t_j)}_s^2,\\
\Norm{T_jt_j}_{s(\Omega\setminus K_{j,m-1})}^2&\leq\Norm{T_jt_j-\pi T_jt_j}_{s(\Omega\setminus K_{j,m-1})}^2+\Norm{\pi T_jt_j}_{s(\Omega\setminus K_{j,m-1})}^2\\
&\leq\frac{1}{\Lambda}\Norm{T_jt_j}_{a(\Omega\setminus K_{j,m-1})}^2+ \Norm{\pi T_jt_j}_{s(\Omega\setminus K_{j,m-1})}^2.
\end{aligned}
$$
Combining with the property of $\pi$  and the definition of $\eta_j^{m-1,m}$ to obtain
$$
\begin{aligned}
\Norm{\eta_j^{m-1,m}T_jt_j}_{a}^2\leq &\Norm*{T_jt_j}_{a(\Omega\setminus K_{j,m-1})}^2+\Norm*{T_jt_j}_{s(\Omega\setminus K_{j,m-1})}^2\\
\leq&\RoundBrackets*{1+\frac{1}{\Lambda}}\Norm*{T_jt_j}_{a(\Omega\setminus K_{j,m-1})}^2+\Norm*{\pi T_jt_j}_{s(\Omega\setminus K_{j,m-1})}^2.
\end{aligned}
$$
Finally,
$$\begin{aligned}
I_3\leq C(\Lambda)\RoundBrackets{\Norm*{(T_j-T_{j,m})t_j}_a^2+\Norm*{(T_j-T_{j,m})t_j}_s^2}^{1/2}\\
\RoundBrackets{\Norm*{T_jt_j}_{a(\Omega\setminus K_{j,m-1})}^2+\Norm*{\pi T_jt_j}_{s(\Omega\setminus K_{j,m-1})}^2}^{1/2}.
\end{aligned}$$
Similarly, we can use the resolution condition for $I_4$,
\begin{equation*}
\begin{aligned}
I_4\leq&c_*k^2H^2\varepsilon^{-2}\RoundBrackets{\Norm*{(T_j-T_{j,m})t_j}_s^2+\Norm{\RoundBrackets{(T_j-T_{j,m})t_j,\eta_j^{m-1,m}T_jt_j}}_s}\\
\leq&\frac{1}{2}\Norm*{(T_j-T_{j,m})t_j}_s^2+\frac{1}{2}\Norm{\RoundBrackets{(T_j-T_{j,m})t_j,\eta_j^{m-1,m}T_jt_j}}_s.\\
\end{aligned}
\end{equation*}
By using the similar tricks in analysing $I_3$,
$$
\begin{aligned}
I_4\leq &\frac{1}{2\Lambda}\Norm{(T_j-T_{j,m})t_j}_a^2+ \frac{1}{2}\Norm{\pi((T_j-T_{j,m})t_j)}_s^2\\+&C(\Lambda)\RoundBrackets{\Norm*{(T_j-T_{j,m})t_j}_a^2+\Norm*{(T_j-T_{j,m})t_j}_s^2}^{1/2}\\&\RoundBrackets{\Norm*{T_jt_j}_{a(\Omega\setminus K_{j,m-1})}^2+\Norm*{\pi T_jt_j}_{s(\Omega\setminus K_{j,m-1})}^2}^{1/2}.    
\end{aligned}
$$
For the remaining terms $I_1$ and $I_2$, we can use Cauchy-Schwarz inequality,
$$I_1\leq\Norm{(T_j-T_{j,m})t_j}_a\Norm{\eta_j^{m-1,m}T_jt_j}_a, I_2\leq\Norm{\pi((T_j-T_{j,m})t_j)}_s\Norm{\pi(\eta_j^{m-1,m}T_jt_j)}_s.$$
We still need to provide a estimate for $\Norm{\pi(\eta_j^{m-1,m}T_jt_j)}_{s}$ by using the property of $\pi$  and the definition of $\eta_j^{m-1,m}$ to obtain
$$
\begin{aligned}
\Norm{\pi(\eta_j^{m-1,m}T_jt_j)}_{s}\leq&\Norm{\eta_j^{m-1,m}T_jt_j}_{s}\leq\Norm*{T_jt_j}_{s(\Omega\setminus K_{j,m-1})}\\
\leq&\frac{1}{\sqrt{\Lambda}}\Norm*{T_jt_j}_{a(\Omega\setminus K_{j,m-1})}+\Norm*{\pi T_jt_j}_{s(\Omega\setminus K_{j,m-1})}.
\end{aligned}
$$
By combining all above together,
$$
\begin{aligned}
&\RoundBrackets*{1-\frac{1}{2\Lambda}}\Norm*{(T_j-T_{j,m})t_j}_{a}^2+\frac{1}{2}\Norm*{\pi(T_j-T_{j,m})t_j}_{s}^2\\
\leq& C(\Lambda)\RoundBrackets{\Norm*{(T_j-T_{j,m})t_j}_{a}^2+\Norm*{\pi(T_j-T_{j,m})t_j}_{s}^2}^{1/2}\\
&\RoundBrackets{\Norm*{T_jt_j}_{a(\Omega\setminus K_{j,m-1})}^2+\Norm*{\pi T_jt_j}_{s(\Omega\setminus K_{j,m-1})}^2}^{1/2}\\
\leq&C(\Lambda\CurlyBrackets*{C_{\mathup{ol}}(m+1)^d(\Norm{z}_{a}^2+\Norm{\pi z}_{s}^2)}^{1/2}\CurlyBrackets*{\sum_{j=1}^N\RoundBrackets*{\mathcal{B}(z_j,z_j)+s(\pi z_j,\pi z_j)}}^{1/2}.
\end{aligned}
$$
We can choose $\Lambda\geq\frac{1}{2}$ to obtain 
$$
\begin{aligned}
&\Norm*{(T_j-T_{j,m})t_j}_{a}^2+\Norm*{\pi(T_j-T_{j,m})t_j}_{s}^2\\
\leq& C(\Lambda)\RoundBrackets{\Norm*{(T_j-T_{j,m})t_j}_{a}^2+\Norm*{\pi(T_j-T_{j,m})t_j}_{s}^2}^{1/2}\\
&\RoundBrackets{\Norm*{T_jt_j}_{a(\Omega\setminus K_{j,m-1})}^2+\Norm*{\pi T_jt_j}_{s(\Omega\setminus K_{j,m-1})}^2}^{1/2}.
\end{aligned}
$$
By using the \cref{lm3}, we can derive the conclusion.
\end{proof}
Before we proceed the next proof, we need a assumption of shape regularity for $K_j\in\mathcal{T}_H$ such that there is a bound $C_{\mathup{ol}}$ and $m>0$,
$$\text{card}\{K\in\mathcal{T}_H: K\subset K_j^m\}\leq C_{\mathup{ol}}m^d.$$
\begin{lemma}\label{lm5}
There exists a constant $0<\beta<1$, independent of $H, m$ and $k$ such that for any $t_j\in V_H$
$$
\Norm*{\sum_{j=1}^N(T_j-T_{j,m})t_j}_a^2+\Norm*{\pi \sum_{j=1}^N(T_j-T_{j,m})t_j}_s^2\leq C^2(\Lambda) C_{\mathup{ol}}\beta^m(m+1)^d\sum_{j=1}^N s(\pi t_j, \pi Tt_j).
$$
\end{lemma}
\begin{proof}
To obtain the global estimate, we set $z :=\sum_{j=1}^N (T_j-T_{j,m})t_j$ and $z_j:=(T_j-T_{j,m})t_j$. Due to supp$(\eta_j^{m,m+1}z)\subset\Omega\setminus K_{j,m}$, supp$(\pi(\eta_j^{m,m+1}z))\subset\Omega\setminus K_{j,m}$, 
supp$(T_{j,m}t_j)\subset\overline{K_{j,m}}$ and supp$(\pi T_{j,m}t_j)\subset\overline{K_{j,m}}$, By using the previous definition of the operators $T_j$ and $T_{j,m}$, after subtracting each other we have 
\begin{equation*}
\mathcal{B}(z_j,\eta_j^{m,m+1}z)+s(\pi z_j,\pi\eta_j^{m,m+1}z)=0.
\end{equation*}
After reformulating the above equation and using the definition of the operator $\mathcal{B}$
\begin{equation}\label{eq:16}
\begin{aligned}
\mathcal{B}\RoundBrackets*{z_j,z}+s\RoundBrackets*{\pi z_j,\pi z}
=&\mathcal{B}(z_j,z-\eta_j^{m,m+1}z)+s(\pi z_j,\pi(z-\eta_j^{m,m+1}z))\\
=&a(z_j,z-\eta_j^{m,m+1}z)+s(\pi z_j,\pi(z-\eta_j^{m,m+1}z))\\
&-\mathrm{i}k(z_j,z-\eta_j^{m,m+1}z)_{\Gamma}-k^2(z_j, z-\eta_j^{m,m+1}z).
\end{aligned}
\end{equation}
For the last four terms of the above, we take the absolute value of them and use the Trace inequality,
\begin{equation*}
\begin{aligned}
\abs{\mathrm{i}k(z_j,z-\eta_j^{m,m+1}z)_{\Gamma}}\leq&
\Norm{z_j}_{1,A,k}\Norm{z-\eta_j^{m,m+1}z}_{1,A,k}\\
\leq&\RoundBrackets*{\Norm{z_j}_a^2+k^2\Norm{z_j}^2}^{1/2}\RoundBrackets{\Norm{z-\eta_j^{m,m+1}z}_a^2+k^2\Norm{z-\eta_j^{m,m+1}z}}^{1/2}.
\end{aligned}
\end{equation*}
by using resolution condition again,
\begin{equation*}
\begin{aligned}
\abs{k^2(z_j,z-\eta_j^{m,m+1}z)}
\leq&c_*k^2H^2\varepsilon^{-2}\RoundBrackets{\Norm{z_j}_s \Norm{z-\eta_j^{m,m+1}z}_s}\leq\frac{1}{2}\Norm{z_j}_s\Norm{z-\eta_j^{m,m+1}z}_s.
\end{aligned}
\end{equation*}
Then taking the absolute value of both hand sides of \cref{eq:16} 
\begin{equation*}
\begin{aligned}
&\abs{\mathcal{B}\RoundBrackets*{z_j,z}+s\RoundBrackets*{\pi z_j,\pi z}}\\
\leq&\Norm{z_j}_a\Norm{z-\eta_j^{m,m+1}z}_a+\Norm{\pi z_j}_s\Norm{\pi(z-\eta_j^{m,m+1}z)}_s+\frac{1}{2}\Norm{z_j}_s\Norm{z-\eta_j^{m,m+1}z}_s\\
&+\RoundBrackets{\Norm{z_j}_a^2+k^2\Norm{z_j}^2}^{1/2}\RoundBrackets{\Norm{z-\eta_j^{m,m+1}z}_a^2+k^2\Norm{z-\eta_j^{m,m+1}z}^2}^{1/2}\\
=&\sum_{i=1}^4I_i.
\end{aligned}
\end{equation*}
We can similarly use the property of $\pi$  and the definition of $\eta_j^{m,m+1}$ to obtain
$$
\begin{aligned}
\Norm{z-\eta_j^{m,m+1}z}_{a}\leq &\Norm{z}_{a(K_{j,m+1})}+\Norm{z}_{s(K_{j,m+1})}\\
\leq&\RoundBrackets{1+\frac{1}{\sqrt{\Lambda}}}\Norm{z}_{a(K_{j,m+1})}+\Norm{\pi z}_{s(K_{j,m+1})}.\\
\Norm{\pi(z-\eta_j^{m,m+1}z)}_{s}\leq&\Norm{z-\eta_j^{m,m+1}z}_{s}\leq\Norm{z}_{s(K_{j,m+1})}\\
\leq&\frac{1}{\sqrt{\Lambda}}\Norm{z}_{a( K_{j,m+1})}+\Norm{\pi z}_{s( K_{j,m+1})}.
\end{aligned}
$$
By using the above, we can obtain the following estimates of $I_1$ and $I_2$, which is 
$$I_1+I_2\leq C(\Lambda)(\Norm{z_j}_a+\Norm{\pi z_j}_s)(\Norm{z}_{a( K_{j,m+1})}+\Norm{\pi z}_{s( K_{j,m+1})}).$$
Due to the following estimates
$$
\begin{aligned}
\Norm{z_j}_s\leq \Norm{z_j-\pi z_j}_s+\Norm{\pi z_j}_s\leq \frac{1}{\sqrt{\Lambda}}\Norm{z_j}_a+\Norm{\pi z_j}_s,
\end{aligned}
$$
then for $I_3$
$$I_3\leq C(\Lambda)(\Norm{z_j}_a+\Norm{\pi z_j}_s)(\Norm{z}_{a(K_{j,m+1})}+\Norm{\pi z}_{s(K_{j,m+1})}).$$
Also,
$$
\begin{aligned}
k^2\Norm{z_j}^2\leq&c_*H^2k^2\varepsilon^{-2} \Norm{z_j}_s^2\leq\frac{1}{2}\Norm{z_j}_s^2\leq \frac{1}{2}\RoundBrackets*{\frac{1}{\Lambda}\Norm{z_j}_a^2+\Norm{\pi z_j}_s^2}.
\end{aligned}
$$
and
$$
\begin{aligned}
k^2\Norm{z-\eta_j^{m,m+1}z}^2&\leq  c_*H^2k^2\varepsilon^{-2} \Norm{z-\eta_j^{m,m+1}z}_s^2\leq \frac{1}{2}\Norm{z}_{s(K_{j,m+1})}^2\\
&\leq \frac{1}{2\Lambda}\Norm{z}_{a( K_{j,m+1})}^2+\frac{1}{2}\Norm{\pi z}_{s( K_{j,m+1})}^2,
\end{aligned}
$$
we can derive the last estimate for $I_4$,
$$I_4\leq C(\Lambda)(\Norm{z_j}_a^2+\Norm{\pi z_j}_s^2)^{1/2}(\Norm{z}_{a( K_{j,m+1})}^2+\Norm{\pi z}_{s( K_{j,m+1})}^2)^{1/2}. $$
Collecting all the estimates of $I_1$ to $I_4$, we hence obtain
\begin{equation*}
\begin{aligned}
\abs{\mathcal{B}\RoundBrackets{z_j,z}+s\RoundBrackets{\pi z_j,\pi z}}\leq C(\Lambda)(\Norm{z}_{a(K_{j,m+1})}^2+\Norm{\pi z}_{s(K_{j,m+1})}^2)^{1/2}(\Norm{z_j}_{a}^2+\Norm{\pi z_j}_{s}^2)^{1/2}.
\end{aligned}
\end{equation*}
Recalling the definition of $C_{\mathup{ol}}$, it is easy to show
\begin{equation*}
\sum_{j=1}^N\RoundBrackets{\Norm{z}_{a(K_{j,m+1})}^2+\Norm{\pi z}_{s(K_{j,m+1})}^2}\leq C_{\mathup{ol}}(m+1)^d\RoundBrackets*{\Norm{z}_{a}^2+\Norm{\pi z}_{s}^2}.
\end{equation*}
Meanwhile, we get
\begin{equation*}
\begin{aligned}
\Norm{z}_a^2+\Norm{\pi z}_{s}^2=&\mathcal{B}\RoundBrackets*{z,z}+s\RoundBrackets*{\pi z,\pi z}+k^2(z,z)+{\mathrm{i}k(z,z)_{\Gamma}}.
\end{aligned}
\end{equation*}
Taking the real part and the absolute value of the above and using the similar tricks from before, 
$$\abs{\text{Re}\RoundBrackets*{\mathrm{i}k(z,z)_{\Gamma}}}=0,\quad \abs{k^2(z,z)}\leq c_*k^2H^2\varepsilon^{-2}\Norm{z}_s^2\leq\frac{1}{2}\RoundBrackets*{\frac{1}{\Lambda}\Norm{z}_a^2+\Norm{\pi z}_s^2}.$$
Then 
\begin{equation*}
\begin{aligned}
\RoundBrackets*{1-\frac{1}{2\Lambda}}\Norm{z}_a^2+\frac{1}{2}\Norm{\pi z}_{s}^2
\leq \abs{\mathcal{B}(z,z)+s\RoundBrackets*{\pi z,\pi z} }.
\end{aligned}
\end{equation*}
We also need to choose $\Lambda>\frac{1}{2}$ to obtain
\begin{equation}\label{eq:17}
\Norm{z}_a^2+\Norm{\pi z}_{s}^2
\leq C(\Lambda)\abs{\mathcal{B}(z,z)+s\RoundBrackets*{\pi z,\pi z} }.
\end{equation}
Finally, we can derive the global estimates
\begin{equation*}
\begin{aligned}
\Norm{z}_a^2+\Norm{\pi z}_{s}^2\leq & C(\Lambda)\abs{\mathcal{B}(z,z)+s\RoundBrackets*{\pi z,\pi z} }\\
\leq & C(\Lambda)\abs{\sum_{j=1}^N \mathcal{B}\RoundBrackets*{z_j,z}+s\RoundBrackets*{\pi z_j,\pi z} }\\
\leq& C(\Lambda) \sum_{j=1}^N\abs{\mathcal{B}\RoundBrackets*{z_j,z}+s\RoundBrackets*{\pi z_j,\pi z} }\\
\leq &C(\Lambda)\CurlyBrackets*{\sum_{j=1}^N\RoundBrackets*{\Norm{z_j}_{a}^2+\Norm{\pi z_j}_{s}^2}}^{1/2}\CurlyBrackets*{\sum_{j=1}^N (\Norm{z}_{a(K_{j,m+1})}^2+\Norm{\pi z}_{s(K_{j,m+1})}^2)}^{1/2}\\
\leq&C(\Lambda)\CurlyBrackets*{C_{\mathup{ol}}(m+1)^d(\Norm{z}_{a}^2+\Norm{\pi z}_{s}^2)}^{1/2}\\&\CurlyBrackets*{C(\Lambda)\beta^m\sum_{j=1}^N(\Norm{T_jt_j}_a^2+\Norm{\pi T_jt_j}_s^2)}^{1/2}\\
\leq&C(\Lambda)\CurlyBrackets*{C_{\mathup{ol}}(m+1)^d(\Norm{z}_{a}^2+\Norm{\pi z}_{s}^2)}^{1/2}\\&\CurlyBrackets*{C(\Lambda)\beta^m\sum_{j=1}^N\abs{\mathcal{B}(T_jt_j,T_jt_j)+s(\pi T_jt_j,\pi T_jt_j)}}^{1/2},\\
\end{aligned}
\end{equation*}
where the last inequality comes from the \cref{lemma2} and \cref{eq:17} by substituting $z$ to $T_jt_j$. By using the previous results
$$\mathcal{B}(T_jt_j,T_jt_j)+s(\pi T_jt_j,\pi T_jt_j)=s(\pi t_j, \pi T_jt_j)\leq\Norm{\pi t_j}_s^2.,$$
we can obtain the desired proof.
\end{proof}
Next, we are going to prove that the global basis functions are indeed localizable. For this purpose, we need to define a bubble function $B(x)$ satisfies the following: for each coarse block $K_j\in\mathcal{T}_{H}$, then 
$$
\begin{cases}
B(x)>0,\quad&\forall x \in\text{int}(K_j),\\
B(x)=0,\quad&\forall x \in\partial K_j.
\end{cases}
$$
We take $B=\Pi_j\eta_j$ where the product is taken over all vertices $j$ on $\partial K_j$. Using the bubble function, we define the following constant 
$$
C_{\pi}=\sup_{K\in\mathcal{T}, \mu\in V_{\mathup{aux}}}\frac{\int_{K}\tilde{\kappa}\mu^2\dx x}{\int_KB\tilde{\kappa}\mu^2\dx x}.
$$
In the following, we are going to prove \cref{lm4.7}, which says that for any $\psi\in V_\mathup{aux}$, we can find a function $\phi$ in the space $V$ such that the $a$-norm of the function $\phi$ is controlled by the $s$-norm of $\psi$ and the support of the function $\phi$ is contained in the support of the function $\psi$.
\begin{lemma}\label{lm4.7} For all $v_\mathup{aux}\in V_\mathup{aux}$, there exists a function $v\in H_0^1(\Omega)$ such that 
$$\pi(v)=v_{\mathup{aux}},\quad\Norm{v}_a\leq C\Norm{v_{\mathup{aux}}}_s^2,\quad \mathup{supp}(v)\subset\mathup{supp}(v_{\mathup{aux}}).
$$
\end{lemma}
\begin{proof}
Consider the space $V_\mathup{aux}(K_j)$ and $\psi\in V_\mathup{aux}(K_j)$  , then we need to find $\phi\in V(K_j)$ and $\lambda\in V_\mathup{aux}$ such that 
$$
\begin{aligned}
\int_{K_j} A\nabla\psi\cdot\nabla\overline{p} \dx x +\int_{K_j}\tilde{\kappa}p\overline{\lambda}\dx x=&0 \quad&&\forall p\in V(K_{j}),\\
\int_{K_j}\tilde{\kappa}\psi\overline{q}\dx x=&\int_{K_j}\tilde{\kappa}\phi\overline{q}\dx x\quad&&\forall q \in V_\mathup{aux}(K_{j}).
\end{aligned}
$$
The above two equations is equivalent to the following minimization problem defined on the coarse block $K_j$:
for a given $\psi\in V_\mathup{aux}^j$ with $\Norm{\psi}=1$, then
\begin{equation}\label{local}
\phi=\text{argmin}\CurlyBrackets*{a(w,w)\,|\, w\in V_0(K_j),s_j(w,\psi)=1, s_j(w,m)=0, \forall m\in V_\mathup{aux}^{\perp}}.
\end{equation}
The well-posedness of the problem \cref{local} is equivalent to the existence of a function $\phi\in V_0(K_j)$ such that 
$$s_j(\phi,\psi)\geq C\Norm{\psi}_s^2,\quad\Norm{\phi}_{a(K_j)}\leq C\Norm{\psi}_{s(K_j)},$$
where $C$ is independent of the meshsize but possibly depends on the problem parameters.
Note that $\psi$ is supported in $K_j$. We let $\phi=B\psi$. By the definition of $s_j$, we have 
$$
s_j(\psi,\phi)=\int_{K_j}\tilde{\kappa}B\abs{\psi}^2\dx x\geq C_{\pi}^{-1}\Norm{\psi}_{s(K_j)}^2.
$$
Since 
$$\nabla(B\psi)=\psi\nabla B+B\nabla\psi, \quad\abs{B}\leq 1,\quad\abs{\nabla B}\leq C_T\sum_j\abs{\nabla\eta_j}^2,
$$
then we have 
$$
\Norm{\phi}_{a(K_j)}^2=\Norm{B\psi}_{a(K_j)}^2\leq C_TC_{\pi}\Norm{\phi}_{a(K_j)}\RoundBrackets*{\Norm{\psi}_{a(K_j)}+\Norm{\psi}_{s(K_j)}}.
$$
Finally, using the spectral problem \cref{eig}, we can obtain
$$
\Norm{\psi}_{a(K_j)}^2\leq \RoundBrackets*{\max_{1\leq i\leq l_j}\lambda_i^j}\Norm{\psi}_{s(K_j)}.
$$
This proves the unique solvability of the minimization problem \cref{local}.
\end{proof}
The following lemma gives the well-posedness of the multiscale problem \cref{cemvar} by proving the inf-sup condition.
\begin{theorem}\label{thm:wellmul}
Under the resolution condition in \cref{lm1} and the oversampling condition $m\geq \abs{\log(k\varepsilon^{-1})}$, the bilinear form $\mathcal{B}$ satisfies the following inf-sup condition: there exists a constant $\gamma_{\mathup{cem}}>0$ depends on $k$ such that
\[
\inf_{u_H\in V_{\mathup{ms}}} \sup_{v_H^*\in V_{\mathup{ms}}^{*}}\frac{\Re\mathcal{B}(u_H,v^*_H)}{\Norm{u_H}_{1,A,k}\Norm{v_H^*}_{1,A,k}}\geq \gamma_{\mathup{cem}}.
\]
\end{theorem}
\begin{proof}
For $u_H\in V_{\mathup{ms}}$, we can find $u_H=T_m\psi$. Next we choose $u_H'=T\psi\in V_{\mathup{glo}}$. Recalling the inf-sup stability on $V_{\mathup{glo}}$ and $V_{\mathup{glo}}^*$, we can find $v_H'=T^*\phi\in V_{\mathup{glo}}^*$ such that $\mathup{Re}\mathcal{B}(u_H',v_H')\geq\gamma \Norm{u_H'}_{1,A,k}\Norm{v_H'}_{1,A,k}$. Similarly, if we denote $v_H=T_m^*\phi$, we have 
Then we have 
$$
\begin{aligned}
\mathcal{B}(u_H,v_H)=&\mathcal{B}(u_H',v_H')+\mathcal{B}(u_H',v_H-v_H')+\mathcal{B}(u_H-u_H',v_H)\\
\mathup{Re}\mathcal{B}(u_H,v_H)=&\mathup{Re}\mathcal{B}(u_H',v_H')+\mathup{Re}\mathcal{B}(u_H',v_H-v_H')+\mathup{Re}\mathcal{B}(u_H-u_H',v_H)\\
\geq&\gamma\Norm{u_H'}_{1,A,k}\Norm{v_H'}_{1,A,k}-\Norm{u_H'}_{1,A,k}\Norm{v_H-v_H'}_{1,A,k}\\
&-\Norm{u_H-u_H'}_{1,A,k}\Norm{v_H}_{1,A,k}\\
\geq&\gamma\Norm{u_H'}_{1,A,k}\Norm{v_H'}_{1,A,k}-\Norm{u_H'}_{1,A,k}\Norm{(T^*-T_m^*)\phi}_{1,A,k}\\
&-\Norm{(T-T_m)\psi}_{1,A,k}\Norm{v_H}_{1,A,k}.
\end{aligned}
$$
The next we need to show that
$$\Norm{(T^*-T^*_m)\phi}_{1,A,k}\leq\Norm{T^*\phi}_{1,A,k}\,\,\text{and}\,\, \Norm{(T-T_m)\psi}_{1,A,k}\leq\Norm{T\psi}_{1,A,k}.$$
Combining \cref{lm5} to obtain 
$$
\begin{aligned}
\Norm{(T_m-T)\psi}_{a}^2+k^2\Norm{(T_m-T)\psi}^2
\leq C(\Lambda)\sqrt{C_{\mathup{ol}}}\beta^{m/2}(m+1)^{d/2}\Norm{\pi\psi}_s^2.
\end{aligned}
$$
where we use the resolution condition
\begin{equation*}
\begin{aligned}
k^2\Norm{(T_m-T)\psi}^2\leq c_*k^2H^2\varepsilon^{-2}\Norm{(T_m-T)\psi}_s^2
\leq\frac{1}{2\Lambda}\Norm{(T_m-T)\psi}_a^2+\Norm{\pi (T_m-T)\psi}_s^2.
\end{aligned}
\end{equation*}
Our remaining task is to $\Norm{\pi\psi}_s$ can be bounded by $\Norm{T\psi}_a$ due to the following
$$\Norm{\pi\psi}_s^2+k^2\Norm{T\psi}^2\leq\Norm{T\psi}_{a}^2+k^2\Norm{T\psi}^2=\Norm{T\psi}_{1,A,k}^2.$$
Take $\psi$ as the test function, we can obtain 
$$
\mathcal{B}(T\psi,\psi)+s(\pi T\psi,\pi\psi)=s(\pi \psi,\pi\psi).
$$
Then
$$
\begin{aligned}
\Norm{\pi\psi}_s^2&\leq\Norm{T\psi}_{1,A,k}\Norm{\psi}_{1,A,k}+\Norm{\pi T\psi}_s\Norm{\pi\psi}_s\\
&\leq(\Norm{T\psi}_a^2+k^2\Norm{T\psi}^2)^{1/2}(\Norm{\psi}_a^2+k^2\Norm{\psi}^2)^{1/2}+\Norm{\pi T\psi}_s\Norm{\pi\psi}_s.
\end{aligned}$$
We again use the resolution condition
\begin{equation*}
\begin{aligned}
k^2\Norm{T\psi}^2&\leq c_*k^2H^2\varepsilon^{-2}\Norm{T\psi}_s^2
\leq\frac{1}{2\Lambda}\Norm{T\psi}_a^2+\Norm{\pi T\psi}_s^2,\\
k^2\Norm{\psi}^2&\leq c_*k^2H^2\varepsilon^{-2}\Norm{\psi}_s^2
\leq\frac{1}{2\Lambda}\Norm{\psi}_a^2+\Norm{\pi\psi}_s^2.
\end{aligned}
\end{equation*}
By using the fact in \cref{lm4.7} to obatin $\Norm{\psi}_a\leq\Norm{\pi\psi}_s$, finally we can obtain
$$
\begin{aligned}
\Norm{\pi\psi}_s^2&\leq\Norm{T\psi}_{1,A,k}\Norm{\psi}_{1,A,k}+\Norm{\pi T\psi}_s\Norm{\pi\psi}_s\\
&\leq(C(\Lambda)\Norm{T\psi}_a^2+\Norm{\pi T\psi}_s^2)^{1/2}C(\Lambda)\Norm{\pi\psi}_s+\Norm{\pi T\psi}_s\Norm{\pi\psi}_s.
\end{aligned}$$
Meanwhile, utilizing Poincare inequaliy, we can obtain
$$\Norm{\pi T\psi}_s\leq\Norm{T\psi}_s\leq C_{\mathup{po}}H^{-1}\Norm{T\psi}_a.$$
Then we obtain the desired estimate
$$
\begin{aligned}
\Norm{(T_m-T)\psi}_{a}+k^2\Norm{(T_m-T)\psi}
\leq C(\Lambda)H^{-1}\sqrt{C_{\mathup{ol}}}\beta^{m/2}(m+1)^{d/2}\Norm{T\psi}_{1,A,k}.
\end{aligned}
$$
Using the same tricks to obtain
$$
\begin{aligned}
\Norm{(T_m^*-T^*)\phi}_{a}+k^2\Norm{(T_m^*-T^*)\phi}
\leq C(\Lambda)H^{-1}\sqrt{C_{\mathup{ol}}}\beta^{m/2}(m+1)^{d/2}\Norm{T^*\phi}_{1,A,k}.
\end{aligned}
$$
Finally we obtain 
$$
\Re\mathcal{B}(u_H,v_H^*)\geq(\gamma-2C(\Lambda)H^{-1}\sqrt{C_{\mathup{ol}}}\beta^{m/2}(m+1)^{d/2})\Norm{u_H}_{1,A,k}\Norm{v_H^*}_{1,A,k}.
$$
\end{proof}

\begin{theorem}\label{thm:main}
Let $u$ be the solution of the problem \cref{eq:ell1} and $u_{\mathup{ms}}$ be the solution of the multiscale problem \cref{cemvar}. Then we have
\[
\Norm{u-u_{\mathup{ms}}}_{a}\leq \frac{1}{\varepsilon_0\sqrt{\Lambda}}\Norm{f}_{s^{-1}}+C(\Lambda)\sqrt{C_{\mathup{ol}}}\beta^{m/2}(m+1)^{d/2}(C_{\mathup{inv}}\Norm{u^{\mathup{glo}}}_{a}+\Norm{\pi u^{\mathup{glo}}}_s).
\]
\end{theorem}
\begin{proof}
Denote $e:=u-u_{\mathup{ms}}$, the triangle inequality gives 
$$
\begin{aligned}
\Norm{e}_{a}^2&\leq\Norm{u-u_{\mathup{glo}}}_{a}^2+\Norm{u_{\mathup{glo}}-u_{\mathup{ms}}}_{a}^2\leq\frac{1}{\varepsilon_0\sqrt{\Lambda}}\Norm{f}_{s^{-1}}+\Norm{u_{\mathup{glo}}-u_{\mathup{ms}}}_{a}^2,
\end{aligned}
$$
where \cref{lm2} is applied to the first term of above. The next we estimate the remaining term $\Norm{u_{\mathup{glo}}-u_{\mathup{ms}}}_{a}$. 
Due to the Cea's Lemma, there exist $\psi_*, \psi_*'$ such that 
$u_{\mathup{glo}}=T\psi_*$, 
$u_{\mathup{ms}}=T_m\psi_*'$, and
$$\|T\psi_*-T_m\psi_*'\|_{1,A,k}\leq \inf_{\psi_{*}^{''}\in V}\|T\psi_*-T_m\psi_{*}^{''}\|_{1,A,k}.$$
Then, choose $\psi_{*}{''}=\psi_*\in V$, then 
$$\Norm{u_{\mathup{glo}}-u_{\mathup{ms}}}_{1,A,k}^2\leq \Norm{(T-T_{m})\psi_*}_{1,A,k}^2=\Norm{(T-T_{m})\psi_*}_{a}^2+k^2\Norm{(T-T_{m})\psi_*}^2.$$
Using the resolution condition
\begin{equation*}
\begin{aligned}
k^2\Norm{(T-T_{m})\psi_*}^2\leq &c_*k^2H^2\varepsilon^{-2}\Norm{(T-T_{m})\psi_*}_s^2\\
\leq&\frac{1}{2{\Lambda}}\Norm{(T-T_{m})\psi_*}_a^2+\Norm{\pi (T-T_{m})\psi_*}_s^2.
\end{aligned}
\end{equation*}
By using the same tricks in proving \cref{lm5} to obtain 
$$
\Norm{(T-T_{m})\psi_*}_{a}^2+\Norm{\pi(T-T_{m})\psi_*}_{s}^2\leq 
C^2(\Lambda)C_{\mathup{ol}}\beta^m(m+1)^d\Norm{\pi \psi_*}_{s}^2.
$$
Then we let $v=\psi_*$ and to obtain
$$\Norm{\pi\psi_*}_s^2=\mathcal{B}(u^{\mathup{glo}},\psi_*)+s(\pi u^{\mathup{glo}},\pi\psi_*).$$
According to \cref{lm4.7}, it is possible to find $\hat{\psi_*}\in V$ such that $\hat{\psi_*}=\pi\psi_*$ and $\Norm{\hat{\psi_*}}_a\leq C_{\mathup{inv}}\Norm{\pi\psi_*}_s$.
Recalling the definition of the operator $\mathcal{B}$
$$\mathcal{B}(u^{\mathup{glo}},\hat{\psi_*})=a(u^{\mathup{glo}},\hat{\psi_*})-k^2(u^{\mathup{glo}},\hat{\psi_*})-\mathrm{i}k(u^{\mathup{glo}},\hat{\psi_*})_{\Gamma}.$$
For the second term on the right hand side, using the resolution condition
\begin{equation*}
\begin{aligned}
\abs{k^2(u^{\mathup{glo}},\hat{\psi_*})}\leq &c_*k^2H^2\varepsilon^{-2}\Norm{u^{\mathup{glo}}}_s\Norm{\hat{\psi_*}}_s\\
\leq&\frac{1}{2}\RoundBrackets*{\frac{1}{\sqrt{\Lambda}}\Norm{u^{\mathup{glo}}}_a+\Norm{\pi u^{\mathup{glo}}}_s}\RoundBrackets*{\frac{1}{\sqrt{\Lambda}}\Norm{\hat{\psi_*}}_a+\Norm{\pi \hat{\psi_*}}_s}.
\end{aligned}
\end{equation*}
Similarly,
\begin{equation*}
\begin{aligned}
\abs{\mathrm{i}k(u^{\mathup{glo}},\hat{\psi_*})_{\Gamma}}
\leq& \Norm{u^{\mathup{glo}}}_{1,A,k}\Norm{\hat{\psi_*}}_{1,A,k}\\\leq&\RoundBrackets{\Norm{u^{\mathup{glo}}}_a^2+k^2\Norm{u^{\mathup{glo}}}^2}^{1/2}\RoundBrackets{\Norm{\hat{\psi_*}}_a^2+k^2\Norm{\hat{\psi_*}}^2}^{1/2}\\
\leq&\RoundBrackets*{\Norm{u^{\mathup{glo}}}_a^2+\frac{1}{2}\RoundBrackets*{\frac{1}{\Lambda}\Norm{u^{\mathup{glo}}}_a^2+\Norm{\pi u^{\mathup{glo}}}_s^2}}^{1/2}\\&\RoundBrackets*{\Norm{\hat{\psi_*}}_a^2+\frac{1}{2}\RoundBrackets*{\frac{1}{\Lambda}\Norm{\hat{\psi_*}}_a^2+\Norm{\pi \hat{\psi_*}}_s^2}}^{1/2}.
\end{aligned}
\end{equation*}
All in all, we can obtain
\begin{equation*}
\begin{aligned}
\Norm{\pi\psi_*}_s^2\leq\Norm{\pi\psi_*}_s\RoundBrackets*{C_{\mathup{inv}}\Norm{u^{\mathup{glo}}}_a+C\Norm{\pi u^{\mathup{glo}}}_s}.
\end{aligned}
\end{equation*}
\end{proof}
Moreover, it is easy to obtain $\|f\|_{s^{-1}}=O(H)$.
From the inf-sup condition of the global problem in \cref{thm:well} and \cref{glo}, we obtain for any $v\in V_{\mathup{glo}}^*\subset H^1$,
\begin{equation}\label{eq:4.9.1}
\sup_{v\in H^1}\frac{(f,v)}{\Norm{v}_{1,A,k}}\geq\gamma(k)\Norm{u^{\mathup{glo}}}_{1,A,k}.\end{equation}
By using Cauchy--Schwartz inequality and the definition of the $k$-weighted norm, we obtain from \cref{eq:4.9.1},
$$\gamma(k)\Norm{f}\geq\Norm{u^{\mathup{glo}}}_{1,A,k}.$$
Due to
$$\Norm{u^{\mathup{glo}}}_{1,A,k}^2=\Norm{u^{\mathup{glo}}}_a^2+k^2\Norm{u^{\mathup{glo}}}^2\leq \gamma^2(k)\Norm{f}^2,$$
we choose $\gamma(k)=k^{-1}$ to obtain
\begin{equation*}
\Norm{u^{\mathup{glo}}}_a+\Norm{\pi u^{\mathup{glo}}}_s\leq\Norm{u^{\mathup{glo}}}_a+H^{-1}\Norm{u^{\mathup{glo}}}\leq O(k^{-1})+O((kH)^{-1}). 
\end{equation*}
If we assume $C_{\mathup{inv}}=O(1)$ and choose suitable $m$ such that $\beta^{m/2}(m+1)^{d/2}=O((kH)^2)$, we will have
$$\Norm{u-u_{\mathup{ms}}}_a=O(kH).$$
\section{Numerical experiments}\label{sec:num}
In this section, we provide several examples of Helmholtz equations in different domains solved by using CEM-GMsFEM method, which demonstrates our established theoretical analysis. We let $u_h$ denotes the reference solution and let $e_h$ denotes the error between the reference solution and numerical solution.
The accuracy will be measured both in the $L^2$ norm and energy norm:
\begin{equation*}\label{error}
\begin{aligned}
e_{L^2}:=&\frac{\Norm{e_h}_{L^2(\Omega)} }{\Norm{u_h}_{L^2(\Omega)}},\quad
e_{a}:=&\frac{\Norm{e_h}_{a(\Omega)} }{\Norm{u_h}_{a(\Omega)}}.
\end{aligned}
\end{equation*}
In order to clearly state the experimental results, we list the notations in the following \cref{tab:1}.

We conduct all numerical experiments on a square domain $\Omega=[0,1]\times[0,1]$. We will calculate reference solutions on a $200\times 200$ mesh with the bilinear Lagrange FEM, therefore the media term $A(x)$ is generated from $200\mathup{px}\times 200\mathup{px}$ figures. For the coarse grid , we choose $H$ to be $1/10, 1/20$ and $1/40$.
For simplicity, we take $\tilde{A}$ as $\tilde{A}|_{K_j}=24H^{-2}A|_{K_j}$
for all numerical experiments as suggested in \cite{Ye2022} and we implement all simulations using the Python libraries Numpy and SciPy.
\begin{table}[ht]
\caption{Simplified description of symbols}
\centering
\begin{tabular}{c c}
\hline
Parameters & Symbols \\
\hline
Number of oversampling layers & $m$ \\
Number of basis functions in every coarse element & $N_{\mathup{bf}}$ \\
Length of every coarse element size & $H$ \\
Length of every fine element size & $h$ \\
\hline
\end{tabular}
\label{tab:1}
\end{table}

\subsection{Model problem 1}
In the first model, we consider about the following Helmholtz equation in  with the boundary  and the source term $f(\vec{x})=0$. 
Then we obtained the following equation from \cref{eq:ell1}
\begin{equation}
\left\{
\begin{aligned}
-\nabla\cdot(A(x_1,x_2)\nabla u)-k^2u=&0 \quad &&\forall(x_1,x_2)\in\Omega, \\
A(x_1,x_2)\nabla u\cdot \bm{n}-\mathrm{i}ku =&g(x_1,x_2) \quad &&\forall(x_1,x_2)\in\partial\Omega,
\end{aligned}
\right.
\end{equation}We firstly consider about the homogeneous coefficients $A(\vec{x})=I $ with wave number $k=2^4$. The boundary data $g$ is chosen in \cref{eq:23} such that the problem admits the plane wave solution $u_*=\text{exp}(\mathrm{i}\vec{k}\cdot\vec{x})$ with $\vec{k}=k(0.6,0.8)$.
\begin{equation}\label{eq:23}
g(x_1,x_2)=
\begin{cases}
-\mathrm{i}1.6k\cdot\text{exp}(\mathrm{i}k\cdot0.8x_2),&\text{on}\quad \{0\}\times(0,1),\\
-\mathrm{i}0.4k\cdot\text{exp}(\mathrm{i}k\cdot(0.6+0.8x_2)),&\text{on}\quad\{1\}\times(0,1),\\
-\mathrm{i}1.8k\cdot\text{exp}(\mathrm{i}k\cdot0.6x_1),&\text{on}\quad(0,1)\times\{0\},\\
-\mathrm{i}0.2k\cdot\text{exp}(\mathrm{i}k\cdot(0.6x_1+0.8)),&\text{on}\quad(0,1)\times\{1\}.
\end{cases}
\end{equation}
\begin{figure}[!ht]
\centering
\includegraphics[width=0.8\textwidth]{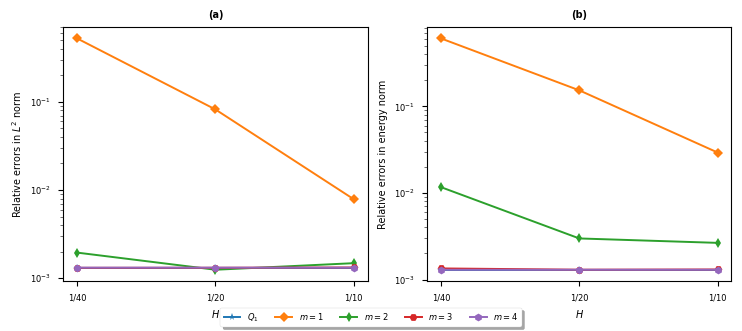}
\caption{Numerical results for the Helmholtz equations with homogeneous coefficients where the relative errors of the proposed method with different numbers of oversampling layers $m$ and the $Q1$ FEM are calculated with respect to the coarse mesh size $H$. Subplots (a) and (b) which the relative errors are measured in the $L_2$ norm and energy norm, respectively.}
\label{fig:sim1}
\end{figure}
\begin{table}[ht]
\centering
\begin{tabular}{ c c c c c }
\hline\hline
$H$  & $N_{\mathup{bf}}$ & $m$ & $\|u_{\mathup{cem}}-u_{h}\|_{L^2(\Omega)}$ & $\|u_{\mathup{cem}}-u_{h}\|_{a(\Omega)}$ \\[0.5ex]
\hline
1/10 &4 & 2 & 1.32e-02 & 1.30e-02 \\
1/20  &4 & 3 & 1.25e-03 & 1.54e-02 \\
1/40  & 4 & 4 & 1.92e-04 & 3.97e-03 \\
\hline
\end{tabular}
\caption{Numerical errors of CEM-GMsFEM using 4 basis functions in each coarse mesh with the oversampling layers equal to $m$. Errors in the $L^2$ norm and energy norm with different $H$.}
\label{tab:3}
\end{table}

For the setting of the proposed multiscale method, we fix $l_j=4$, indicating that we calculate the first four eigenfunctions in \cref{eig} and construct four multiscale bases for each coarse element, while we vary the oversampling layers m from $1$ to $4$. With regard to the relative error norms, we choose $u_h$ to be the exact solution of the Helmholtz equation with the homogeneous coeficients and $u_{\mathup{cem}}$ is approximated by the CEM-GMsFEM method and $e_h$ represents the difference between the $u_h$ and $u_{\mathup{cem}}$. In order to show the efficient of our multiscale method, we also shows the relative error between the exact solution and the approximated solution which is obtained from the $Q1$ FEM. We can observe from subplots (a) to (b) in \cref{fig:sim1} that the convergence of the $Q1$ FEM manifests a linear pattern with respect to $H$ in the logarithmic scale, consistent with the theoretical expectation. We can also see that the number of oversampling layers $m$ has a significant impact on the accuracy of the proposed method and the results of the relative errors with respect to different mesh size and oversampling layers are shown in \cref{tab:3}. However, for the same $m$, the error decaying with respect to $H$ does not always hold ($m=1,2)$, as depicted in subplots (a) and (b). Although for $m = 3,4$, the proposed method exhibits higher accuracy than the $Q1$ FEM, the computational cost is significantly same or higher due to the sophisticated process of constructing multiscale bases. Therefore, the proposed method is more suitable for scenarios involving intricate coefficient profiles.
\subsection{Model problem 2}
\begin{figure}[!ht]
\centering
\begin{subfigure}{0.4\textwidth}
\includegraphics[width=\linewidth]{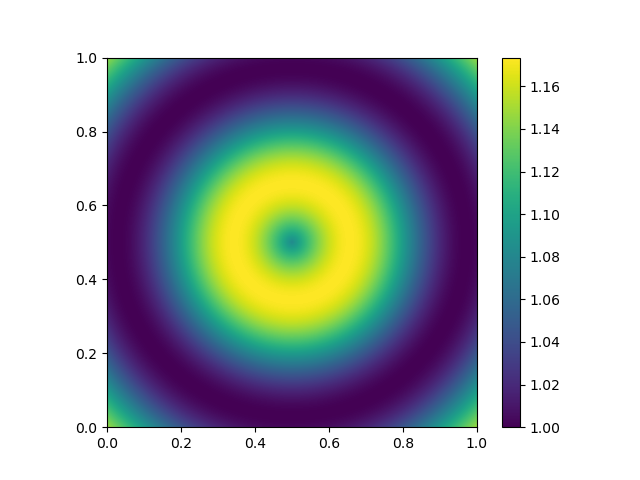} 
\caption{The medium configuration;}
\label{fig:medium1}
\end{subfigure}
\begin{subfigure}{0.4\textwidth}
\includegraphics[width=\linewidth]{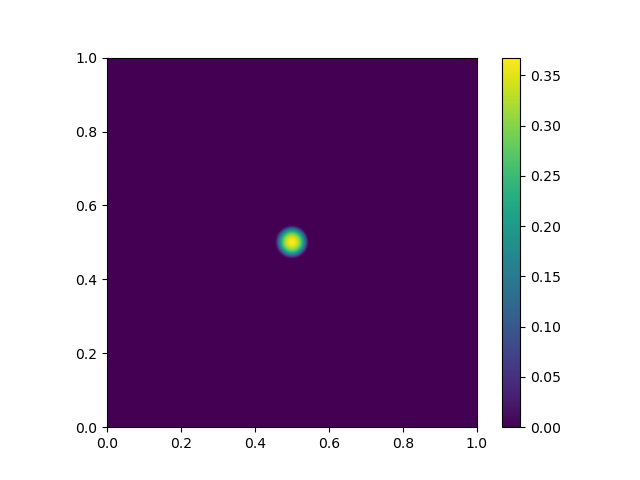}
\caption{The source term $f$.}
\label{fig:source1}
\end{subfigure}
\caption{Model Problem 2 }
\label{md2}
\end{figure}
In the second model, we consider about the following Helmholtz equation with pointwise isotropic coefficients $A(x_1,x_2)$ displayed in \cref{fig:medium1} and source term $f$ is chosen in \cref{eq:22} and it is shown in \cref{fig:source1}.
\begin{equation}\label{eq:22}
f(x_1,x_2)=
\begin{cases}
\exp\RoundBrackets*{-\frac{1}{1-400(x_1^2+x_2^2)}}&\mathup{for}\, \sqrt{x_1^2+x_2^2}<1/20,\\
0&\mathup{else.}
\end{cases}    
\end{equation}
Then we obtained the following equation from \cref{eq:ell1}
\begin{equation}
\left\{
\begin{aligned}
-\nabla\cdot(A(x_1,x_2)\nabla u)-k^2u=&f(x_1,x_2) \quad &&\forall(x_1,x_2)\in\Omega, \\
A(x_1,x_2)\nabla u\cdot \bm{n}-\mathrm{i}ku =&0 \quad &&\forall(x_1,x_2)\in\partial\Omega.
\end{aligned}
\right.
\end{equation}
\begin{table}[!ht]
\centering
\begin{tabular}{ c c c c c }
\hline\hline
$H$ & $N_{\mathup{bf}}$ & $m$ & $\|u_{\mathup{cem}}-u_{h}\|_{L^2(\Omega)}$ & $\|u_{\mathup{cem}}-u_{h}\|_{a(\Omega)}$ \\[0.5ex]
\hline
1/10 & 4 & 2 & 1.80e-02 & 8.24e-02 \\
1/20 & 4 & 3 & 1.00e-03 & 9.27e-03 \\
1/40 & 4 & 4 & 1.03e-04 & 1.87e-03 \\
\hline
\end{tabular}
\caption{Numerical errors of CEM-GMsFEM using 4 basis functions in each coarse mesh with the oversampling layers equal to $m$. Errors in the $L^2$ norm and energy norm with different $H$.}
\label{tab:4}
\end{table}
\begin{figure}[!ht]
\centering
\includegraphics[width=0.8\textwidth]{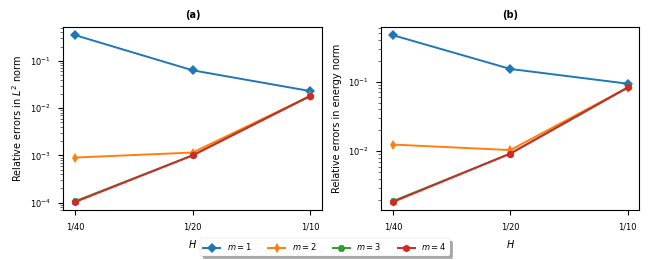}
\caption{Subplots (a) and (b) show the relative errors of the proposed method for the pointwise isotropic coefficients with different numbers of oversampling layers $m$ with respect to the coarse mesh size $H$, but measured in different norms.}
\label{fig:ep2}
\end{figure}

\begin{figure}[!ht]
\centering
\includegraphics[width=0.8\textwidth]{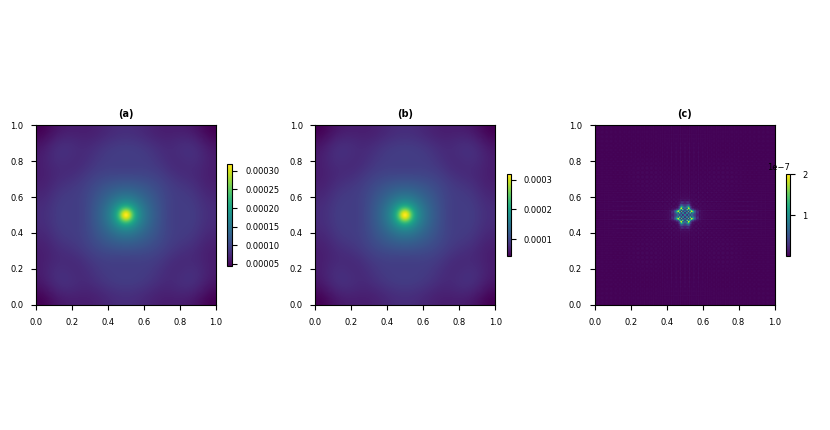}
\caption{Solutions by using $H=1/40$ and oversampling layer equals to 3. Subplot (a), (b),(c) represent reference solution , CEM-GMsFEM solution and Difference of these two solutions, respectively .}
\label{fig:ep2.1}
\end{figure}
The numerical results of the $Q1$ FEM and the proposed method with $m \in\{1, 2, 3\}$ are presented in \cref{fig:ep2}  where Subplots (a) and (b) in \cref{fig:ep2} share the same setting (corresponding to the \cref{fig:sim1}) and the relative errors are measured in the different norms. We choose $u_h$ to be the $Q1$ FEM solution and $u_{\mathup{cem}}$ is approximated by the CEM-GMsFEM method and $e_h$ represents the difference between the $u_h$ and $u_{\mathup{cem}}$, which are shown in \cref{fig:ep2.1}.  The classical CEM-GMsFEM \cite{Yalchin2013} is proven to be effective in handling long and high-contrast channels, and the proposed method inherits this advantage. By setting $m = 3, 4$, the proposed method can achieve a relative error of 10\% in the energy norm. Typically, the relative errors in the $L_2$ norm are significantly smaller by an order of magnitude than those in the energy norm, and the proposed method can achieve a relative error of 0.015 in the $L_2$ norm for $m = 3$. Furthermore, comparing subplots (a) and (b) reveals similar convergence behavior, indicating that the scale of the models does not affect the accuracy of the proposed method. We also list results of the relative errors with respect to different mesh size and oversampling layers are shown in \cref{tab:4}. Note that we solve the Helmholtz equation in coarser mesh, which means we can solve the issues of the pollution effect.
\subsection{Model problem 3}
\begin{figure}[!ht]
\centering
\begin{subfigure}{0.4\textwidth}
\includegraphics[width=\linewidth]{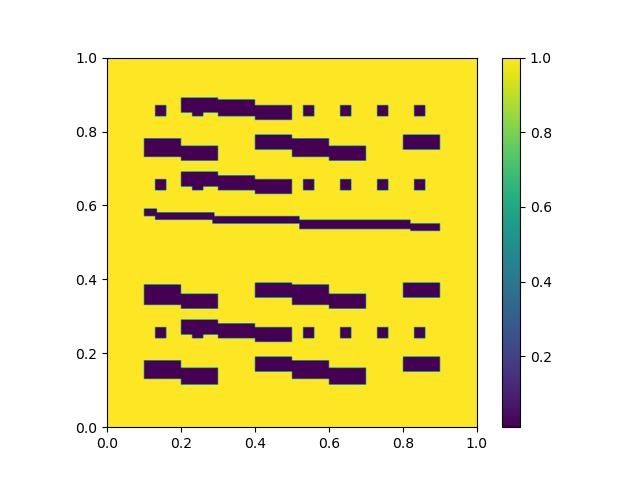} 
\caption{The medium configuration;}
\label{fig:medium2}
\end{subfigure}
\begin{subfigure}{0.4\textwidth}
\includegraphics[width=\linewidth]{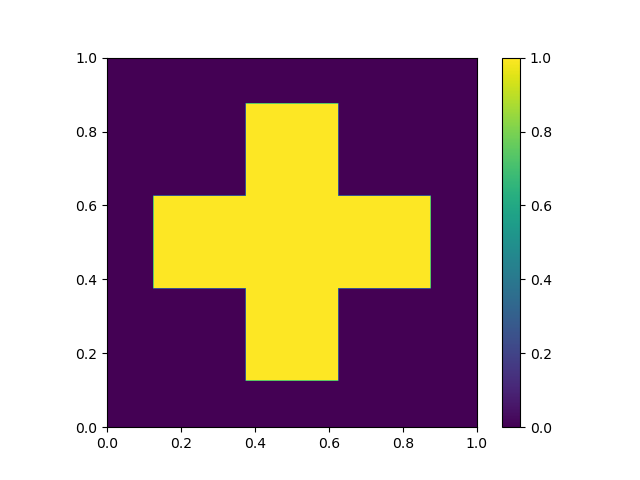}
\caption{The source term $f$.}
\label{fig:source2}
\end{subfigure}
\caption{Model Problem 3 }
\label{md3}
\end{figure}
\begin{table}[!ht]
\centering
\begin{tabular}{ c c c c c }
\hline\hline
$H$ & $N_{\mathup{bf}}$ & $m$ & $\|u_{\mathup{cem}}-u_{h}\|_{L^2(\Omega)}$ & $\|u_{\mathup{cem}}-u_{h}\|_{a(\Omega)}$ \\[0.5ex]
\hline
1/10 & 4 & 2 & 0.9935 & 0.9981 \\
1/20 & 4 & 3 & 0.8260 & 0.8361 \\
1/40 & 4 & 4 & 0.0126 & 0.0213 \\
\hline
\end{tabular}
\caption{Numerical errors of CEM-GMsFEM using 4 basis functions in each coarse mesh with the oversampling layers equal to $m$. Errors in the $L^2$ norm and energy norm with different $H$.}
\label{tab:5}
\end{table}
\begin{figure}[!ht]
\centering
\includegraphics[width=0.8\textwidth]{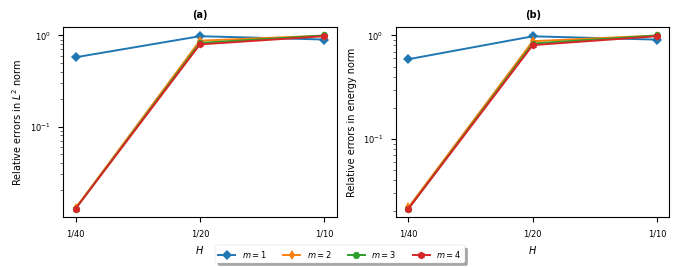}
\caption{Subplots (a) and (b) show the relative errors of the proposed method for the high-contrast coefficients with different numbers of oversampling layers $m$ w.r.t. the coarse mesh size $H$, but measured in different norms.}
\label{fig:ep3}
\end{figure}
\begin{figure}[!ht]
\centering
\includegraphics[width=0.8\textwidth]{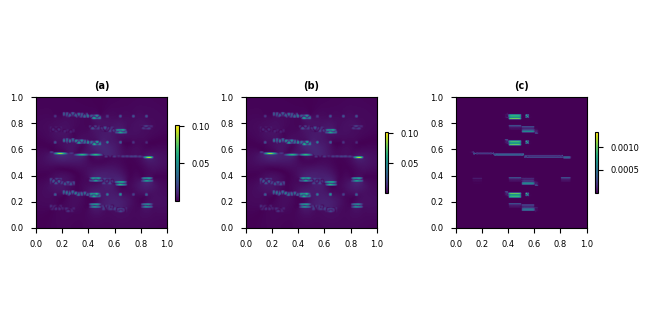}
\caption{Solutions by using $H=1/40$ and oversampling layer equals to 3. Subplot (a), (b),(c) represent reference solution , CEM-GMsFEM solution and difference of above two solutions, respectively.}
\label{fig:ep3.1}
\end{figure}
In the third model, we make use of high-contrast coefficients $A(x_1,x_2)$ displayed in \cref{fig:medium2} with a $10^{-3}$ contrast ratio as the show case for the Helmholtz equation, and source term $f$ is chosen as a piecewise constant function which is shown in \cref{fig:source2}. The numerical results of the proposed multiscale method with $m \in\{1, 2, 3, 4\}$ are presented in \cref{fig:ep3} where subplots (a) and (b) in \cref{fig:ep3} also share the same setting with the mesh seize and wavenumber and the relative errors are measured in the different norms as well. We choose $u_h$ to be the reference solution approximated by the $Q1$ FEM method due to the lack of the classical solution and $e_h$ represents the difference between the multiscale solution $u_{\mathup{cem}}$ and  $Q1$ FEM solution $u_h$, which are shown in \cref{fig:ep3.1}. The proposed method inherits this advantage of CEM-GMsFEM which is  proven to be effective in handling long and high-contrast channels. Specifically, we outlines the errors in both $L^2$ and energy norms corresponding to different chosen $m$ in \cref{tab:5}. By setting $m = 3,4$, the proposed method can achieve a relative error of almost 2\% in the energy norm. Typically, the relative errors in the $L_2$ norm are significantly smaller by an order of magnitude than those in the energy norm, and the proposed method can achieve a relative error of 1.2\% in the $L_2$ norm for $m = 2,3,4$. Furthermore, comparing subplots (a) and (b) in \cref{fig:ep3} reveals similar convergence behavior, indicating that the scale of the models does not affect the accuracy of the proposed method.

\section{Conclusion}\label{con}
The paper introduces the CEM-GMsFEM method for solving Helmholtz equations in heterogeneous medium, provides the proof of convergence for this method, and presents numerical experiments to support its effectiveness. The work contributes to the field by offering a novel approach to tackle heterogeneous Helmholtz problems without restrictive assumptions on the coefficient structures.

In future work, we plan to apply the CEM-GMsFEM method to solve Helmholtz equations in other types of domains or in more complex geometries, such as perforated domains. Perforated domains are characterized by having small holes or voids within the domain, and their inclusion introduces additional challenges in accurately capturing the solution behavior. Additionally, we are also interested in exploring more efficient techniques to solve Helmholtz equations with much higher frequencies. These future directions would contribute to advancing the field of multiscale modeling and simulation of Helmholtz equations in diverse and challenging scenarios.
\section{Acknowledgment}
The research of Eric Chung is partially supported by the Hong Kong RGC General Research Fund (Projects: 14305423 and 14305222).

\bibliographystyle{elsarticle-num-names} 
\bibliography{Helmholtz}

\end{document}